\documentclass[11pt, leqno]{article}

\usepackage{euscript,amstext,amsmath,amsthm,a4,amssymb}

\newcommand{\comment}[1]{}

\newtheorem{thm}{Theorem}[section]

\newtheorem{prop}[thm]{Proposition}

\theoremstyle{remark}
\newtheorem{Rem}[thm]{Remark}
\theoremstyle{definition}
\newtheorem{defn}[thm]{Definition}

\title{Equivariant symmetric bilinear torsions \footnote{This work was
partially supported by the Qiushi Foundation.}}
\author{Guangxiang Su \footnote{Chern Institute of Mathematics \& LPMC, Nankai University,
Tianjin 300071, P.R. China. (sugx@mail.nankai.edu.cn)}}
\date{}

\begin{document}

\maketitle
\begin{abstract}
We extend the main result in the previous paper of Zhang and the
author relating the Milnor-Turaev torsion with the complex valued
analytic torsion to the equivariant case.
\end{abstract}

\renewcommand{\theequation}{\thesection.\arabic{equation}}
\setcounter{equation}{0}

\section{Introduction} \label{s0}

Let $F$ be a unitary flat vector bundle  on  a closed Riemannian
manifold $X$. In \cite{RS}, Ray and Singer defined an analytic
torsion associated to $(X,F)$ and proved that it does not depend on
the Riemannian  metric on $X$. Moreover, they conjectured that this
analytic torsion coincides with the classical  Reidemeister torsion
defined using a triangulation on $X$ (cf. \cite{Mi}). This
conjecture was later proved in the celebrated papers of Cheeger
\cite{C} and M\"uller \cite{Mu1}. M\"uller generalized this result
in \cite{Mu2} to the case where $F$ is a unimodular flat vector
bundle on $X$. In \cite{BZ1}, inspired by the considerations of
Quillen \cite{Q}, Bismut and
 Zhang reformulated the  above Cheeger-M\"uller theorem as an equality between the Reidemeister and
 Ray-Singer metrics defined on the determinant of
  cohomology, and proved an extension of
 it   to the case of general  flat vector bundles over $X$.
 The method used in \cite{BZ1} is different from those of Cheeger
 and M\"uller in that it makes use of a deformation by Morse functions introduced by Witten
 \cite{W}   on the de Rham complex.

 On the other hand,   Turaev generalizes the concept of
 Reidemeister torsion to a complex valued invariant whose
 absolute value provides the original Reidemeister torsion, with the help of the so-called Euler structure (cf.
 \cite{T}, \cite{FT}). It is natural to ask whether there exists
 an analytic interpretation of this Turaev torsion.

 Recently, Burghelea and Haller \cite{BH1, BH2}, following a suggestion of M\"uller,  define a
 generalized
 analytic torsion associated to a nondegenerate  symmetric bilinear form on  a
 flat vector bundle over a closed manifold and
 make an explicit conjecture between this generalized analytic
 torsion and the Turaev torsion. Later this conjecture was proved by
 Su and Zhang \cite{SZ}. Also Burghelea and Haller \cite{BH3}, up to a sign,
 proved this conjecture for odd dimensional manifolds, and comments
 were made how to derive the conjecture in full generality in their
 paper.

 In this paper, we will extend the main result in \cite{SZ} to
 the equivariant case, which is closer in spirit to the approach developed by Bismut-Zhang in \cite{BZ2}.

  The rest of this paper is organized as follows. In Section 2, we
  construct the equivariant symmetric bilinear torsions associated with
  equivariant nondegenerate symmetric bilinear forms on a flat vector  bundle. In Section 3,
  we state the main result of this paper. In Section 4, we provide a
  proof of the main result. Section 5 is
  devoted to the proofs of the intermediary results stated in
  Section 4.

  Since we will make substantial use of the results in \cite{BZ1, BZ2, SZ},
  we will refer to \cite{BZ1, BZ2, SZ} for related definitions and
  notations directly  when there will be no confusion.

\section{Equivariant symmetric bilinear torsions associated to the de Rham and Thom-Smale complexes} \label{s2}
\setcounter{equation}{0}

 In this section, for a $G$-invariant nondegenerate
bilinear symmetric form on a complex flat vector bundle over an
oriented closed manifold, we define  two naturally associated
equivariant symmetric bilinear forms on the equivariant determinant
of the cohomology $H^*(M,F)$  with coefficient $F$. One constructed
in a combinatorial way through the equivariant Thom-Smale complex
associated to a equivariant Morse function, and the other one
constructed in an analytic way through the equivariant de Rham
complex.

\subsection{\normalsize Equivariant symmetric bilinear torsion of a finite dimensional complex}
\label{s2.1}

Let $(C,\partial)$ be a  finite cochain  complex
\begin{align}\label{2.1}
\left(C,\partial\right): 0\longrightarrow
C^0\stackrel{\partial_0}{\longrightarrow}
C^1\stackrel{\partial_1}{\longrightarrow}\cdots\stackrel{\partial_{n-1}}{\longrightarrow}C^n\longrightarrow
0,
\end{align}
where each $C ^i$, $0\leq i\leq n$, is a finite dimensional complex
vector space.

Let each $C^i$, $0\leq i\leq n$, admit a nondegenerate symmetric
bilinear form $b_i$. We equip $C$ with the nondegenerate symmetric
bilinear form $b_{C}=\bigoplus_{i=0}^{n}b_i$.

Let $G$ be a compact group. Let $\rho: G\to {\rm End}(C)$ be a
representation of $G$, with values in the chain homomorphisms of $C$
which preserve the bilinear form $b_C$. In particular, if $g\in G$,
$\rho(g)$ preserves the $C^{i}$'s.

Let $\widehat{G}$ be the set of equivalence classes of complex
irreducible representations of $G$. An element of $\widehat{G}$ is
specified by a complex finite dimensional vector space $W$ together
with an irreducible representation $\rho_W: G\to {\rm End}(W)$.

For $W\in\widehat{G}$, set
\begin{align}\label{2.2}
C_W^i={\rm Hom}_G(W,C^i)\otimes W,
\end{align}
\begin{align}\label{2.3}
C_W={\rm Hom}_G(W,C)\otimes W.
\end{align}
Let $\partial_W$ be the map induced by $\partial$ on $C_W$. Then

\begin{align}\label{2.4}
\left(C_W,\partial_W\right): 0\longrightarrow
C^0_W\stackrel{\partial_{0,W}}{\longrightarrow}
C^1_W\stackrel{\partial_{1,W}}{\longrightarrow}\cdots\stackrel{\partial_{n-1,W}}{\longrightarrow}C^n_W\longrightarrow
0
\end{align}
is a chain complex. Thus we obtain the isotypical decomposition,
\begin{align}\label{2.5}
(C,\partial)=\bigoplus_{W\in\widehat{G}}(C_W,\partial_W),
\end{align}
and the decomposition (\ref{2.5}) is orthogonal.

If $E$ is a complex finite dimensional representation space for $G$,
let $\chi(E)$ be the character of the representation. Put
$$\chi(C)=\sum_{i=0}^{n}(-1)^i\chi(C^i),$$
$$e(C)=\sum_{i=0}^{n}(-1)^i{\rm dim}C^i,$$
\begin{align}\label{2.6}
e(C_W)=\sum_{i=0}^{n}(-1)^i{\rm dim}(C_W^i).
\end{align}
By (\ref{2.5}), we get
\begin{align}\label{2.7}
\chi(C)=\sum_{W\in\widehat{G}}e(C_W){\chi(W)\over{\rm rk}(W)}.
\end{align}

If $\lambda$ is a complex line, let $\lambda^{-1}$ be the dual line.
If $E$ is a finite dimensional complex vector space, set
\begin{align}\label{2.8}
{\rm det}E=\Lambda^{\rm max}(E).
\end{align}

Put
$${\rm det}C=\bigotimes_{i=0}^{n}\left({\rm det}C^i\right)^{(-1)^i},$$
\begin{align}\label{2.9}
{\rm det}C_W=\bigotimes_{i=0}^{n}\left({\rm
det}C_W^i\right)^{(-1)^i}.
\end{align}
By (\ref{2.5}), we obtain
\begin{align}\label{2.10}
{\rm det}C=\bigotimes_{W\in\widehat{G}}{\rm det}C_W.
\end{align}

For $0\leq i\leq n$, $C_W^i$ is a vector subspace of $C^i$. Let
$b_{C_W^i}$ be the induced symmetric bilinear form on $C_W^i$. let
$b_{{\rm det}C_W^i}$ be the symmetric bilinear form on ${\rm
det}C_W^i$ induced by $b_{C_W^i}$, and let $b_{({\rm
det}C_W^i)^{-1}}$ be the dual symmetric bilinear form on $({\rm
det}C_W^i)^{-1}$. Also we have symmetric bilinear forms $b_{{\rm
det}C_W}$ on ${\rm det}C_W$ and $b_{{\rm det}C}$ on ${\rm det}C$.

Put
\begin{align}\label{2.11}
{\rm det}(C,G)=\bigoplus_{W\in\widehat{G}}{\rm det}C_W.
\end{align}

\begin{defn}\label{t2.1}
We introduce the formal product
\begin{align}\label{2.12}
b_{{\rm det}(C,G)}=\prod_{W\in\widehat{G}}\left(b_{{\rm
det}C_W}\right)^{\chi(W)\over{\rm rk}(W)}.
\end{align}
\end{defn}

For $W\in\widehat{G}$, let $x_W$, $y_W \in {\rm det}C_W$, $x_W\neq
0$, $y_W\neq 0$. Set $x=\oplus_{W\in\widehat{G}}x_W$,
$y=\oplus_{W\in\widehat{G}}y_W\in{\rm det}(C,G)$. Then by
definition,
\begin{align}\label{2.13}
b_{{\rm det}(C,G)}(x,y)=\prod_{W\in\widehat{G}}\left(b_{{\rm
det}C_W}(x_W,y_W)\right)^{\chi(W)\over{\rm rk}(W)}.
\end{align}
Tautologically, (\ref{2.13}) is an identity of characters on $G$. In
particular
\begin{align}\label{2.14}
b_{{\rm det}(C,G)}(x,y)(1)=\prod_{W\in\widehat{G}}b_{{\rm
det}C_W}(x_W,y_W).
\end{align}
In fact (\ref{2.14}) just says that
\begin{align}\label{2.15}
b_{{\rm det}(C,G)}(1)=b_{{\rm det}C}.
\end{align}

Of course, using the orthogonality of the $\chi_{W}$'s, knowing the
formal product $b_{{\rm det}(C,G)}$ is equivalent to knowing the
symmetric bilinear forms $b_{{\rm det}C_W}$.

Clearly
$$H(C_W,\partial_W)={\rm Hom}_{G}(W,H(C,\partial))\otimes W,$$
\begin{align}\label{2.16}
H(C,\partial)=\bigoplus_{W\in\widehat{G}}H(C_W,\partial_W).
\end{align}
For $W\in\widehat{G}$, we define ${\rm det}H(C_W,\partial_W)$ as in
(\ref{2.9}). Set
\begin{align}\label{2.17}
{\rm det}(H(C,\partial),G)=\bigoplus_{W\in\widehat{G}}{\rm
det}H(C_W,\partial_W).
\end{align}

For $W\in\widehat{G}$, there is a canonical isomorphism (cf.
\cite{KM} and \cite[Section 1a)]{BGS})
\begin{align}\label{2.18}
{\rm det}C_W\simeq{\rm det}H(C_W,\partial_W).
\end{align}
From (\ref{2.18}), we get
\begin{align}\label{2.19}
{\rm det}(C,G)\simeq{\rm det}(H(C,\partial),G).
\end{align}

Let $b_{{\rm det}H(C_W,\partial_W)}$ be the symmetric bilinear form
on ${\rm det}H(C_W,\partial_W)$ corresponding to $b_{{\rm det}C_W}$
via the canonical isomorphism (\ref{2.18}).

\begin{defn}\label{t2.2}
we introduce the formal product
\begin{align}\label{2.20}
b_{{\rm det}(H(C,\partial),G)}=\prod_{W\in\widehat{G}}\left(b_{{\rm
det}H(C_W,\partial_W)}\right)^{\chi(W)\over{{\rm rk}W}}.
\end{align}
\end{defn}
Tautologically, under the identification (\ref{2.19}),
\begin{align}\label{2.21}
b_{{\rm det}(C,G)}=b_{{\rm det}H((C,\partial),G)}.
\end{align}

By an abuse of notation, we will call the formal product $b_{{\rm
det}(C,G)}$ a symmetric bilinear form on ${\rm det}(C,G)$.

\subsection{\normalsize  The Thom-Smale complex of a gradient field}
\label{s2.2}

Let $M$ be a closed smooth manifold, with ${\rm dim}M=n$. For
simplicity, we make the assumption that $M$ is oriented.

Let $(F,\nabla^F)$ be a complex flat vector bundle over $M$ carrying
the flat connection $\nabla^F$.  We make the assumption that $F$
carries a nondegenerate symmetric bilinear form $b^F$.

Let $(F^*,\nabla^{F^*})$ be the dual complex flat vector bundle of
$(F,\nabla^{F})$  carrying   the dual  flat connection
$\nabla^{F^*}$. 

Let $f:M\rightarrow {\bf R}$ be a Morse function. Let $g^{TM}$ be a
Riemannian metric on $TM$ such that the corresponding gradient
vector field $-X=-\nabla f\in \Gamma(TM)$ satisfies the Smale
transversality conditions (cf. \cite{Sm}), that is, the unstable
cells (of $-X$) intersect transversally with the stable cells.

Set
\begin{align}\label{2.22}
B=\{ x\in M; X(x)=0\} .
\end{align}

For any $ {x}\in {B}$, let $W^u( {x})$ (resp. $W^s( {x})$) denote
the unstable (resp. stable) cell at $ {x}$, with respect to $- {X}$.
We also choose an orientation $O_{ {x}}^-$ (resp. $O_{ {x}}^+$) on
$W^u( {x})$ (resp. $W^s( {x})$).

Let  $ {x}$, $ {y}\in {B}$ satisfy the Morse index relation ${\rm
ind}( {y})={\rm ind}( {x})-1$, then $\Gamma( {x}, {y})=W^u( {x})\cap
W^s( {y})$ consists of a finite number of integral curves $\gamma$
of $- {X}$. Moreover, for each $\gamma\in\Gamma( {x}, {y})$, by
using the orientations chosen above, on can define a number
$n_\gamma( {x},{y})=\pm 1$ as in \cite[(1.28)]{BZ1}.

If $ {x}\in {B}$, let $[W^u( {x})]$ be the complex line generated by
$W^u( {x})$. Set
\begin{align}\label{2.23}
C_*(W^u, {F}^*)=\bigoplus_{ {x}\in {B}}[W^u( {x})]\otimes
 {F}^*_{ {x}},
\end{align}
\begin{align}\label{2.24}
C_i(W^u, {F}^*)=\bigoplus_{ {x}\in {B},\ {\rm ind}( {x})=i}[W^u(
{x})]\otimes  {F}^*_{ {x}}.
\end{align}
If $ {x}\in {B}$, the flat vector bundle $ {F}^*$ is canonically
trivialized on $W^u( {x})$. In particular, if $ {x}$, $ {y}\in {B}$
satisfy  ${\rm ind}( {y})={\rm ind}( {x})-1$, and if
$\gamma\in\Gamma( {x}, {y})$, $f^*\in {F}^*_{ {x}}$, let
$\tau_\gamma(f^*)$ be the parallel transport of $f^*\in {F}^*_{
{x}}$ into $ {F}^*_{ {y}}$ along $\gamma$ with respect to the flat
connection $\nabla^{ {F}^*}$.

Clearly, for any $ {x}\in {B}$, there is only a finite number of $
{y}\in {B}$, satisfying together that ${\rm ind}( {y})={\rm ind}(
{x})-1$ and $\Gamma( {x}, {y})\neq \emptyset$.

If $ {x}\in {B}$, $f^*\in {F}^*_{ {x}}$, set
\begin{align}\label{2.25}
\partial (W^u( {x})\otimes
f^*)=\sum_{ {y}\in {B},\ {\rm ind}( {y})={\rm ind}(x)-1}
\sum_{\gamma\in\Gamma( {x}, {y})}n_\gamma( {x}, {y}) W^u(
{y})\otimes \tau_\gamma(f^*).
\end{align}
Then $\partial$ maps $C_i(W^u, {F}^*)$ into $C_{i-1}(W^u, {F}^*)$.
Moreover, one has
\begin{align}\label{2.26}
\partial^2=0.
\end{align}
That is, $(C_*(W^u, {F}^*),\partial)$ forms a chain complex. We call
it the  Thom-Smale complex associated to $( {M}, F, -X)$.

If $ {x}\in {B}$, let $[W^u( {x})]^*$ be the dual line to $W^u(
{x})$. Let $(C^*(W^u, {F}), {\partial})$ be the complex which is
dual to $(C_*(W^u, {F}^*),\partial)$. For $0\leq i\leq n$, one has
\begin{align}\label{2.27}
C^i(W^u, {F})=\bigoplus_{ {x}\in {B},\ {\rm ind}( {x})=i}[W^u(
{x})]^*\otimes
 {F}_{ {x}}.
\end{align}

Let $G$ be a compact group acting on $M$ by smooth diffeomorphisms.
we assume that the action of $G$ lifts to $F$ and preserves the flat
connection of $F$. Then $G$ acts naturally on $H^*(M,F)$. We assume
that $f$ and $g^{TM}$ are $G$-invariant. Then $-X=-\nabla f$ is also
$G$-invariant. We assume that it verifies the smale transversality
conditions.

Clearly $B$ is $G$-invariant. Also if $x\in B$, $g\in G$,
$$g\left(W^{u}(x)\right)=\epsilon_{g}(x)W^{u}(gx),$$
where $\epsilon_{g}(x)=+1$ if $g(W^{u}(x))$ has the same orientation
as $W^{u}(gx)$, $\epsilon_{g}=-1$ if not. Clearly $g$ acts as a
chain homomorphism on $(C_*(W^{u},F^*),\partial)$. The corresponding
dual action of $g$ on $(C^*(W^u,F),\partial)$ is such that
$$g\left(W^u(x)^*\right)=\epsilon_{g}(x)W^u(gx)^*.$$
Then $g$ acts as a chain homomorphism on $(C^*(W^u,F),\partial)$.
Therefore $g$ acts on $H^*(C^*(W^u,F),\partial)$.

\subsection{\normalsize  Equivariant Milnor symmetric bilinear torsion}
\label{s2.3}

For $x\in B$, let $b^{F_x}$ be a nondegenerate symmetric bilinear
form on $F_x$. We assume that the $b^{F_x}$'s are $G$-invariant,
i.e. for $g\in G$, $x\in B$
\begin{align}\label{2.29}
g\left(b^{F_x}\right)=b^{F_{g(x)}}.
\end{align}

The symmetric bilinear forms $b^{F_x}$'s determine a $G$-invariant
symmetric bilinear form on $C^*(W^u,F)=\bigoplus_{x\in
B}[W^u(x)]^*\otimes F_x$, such that the various $[W^u(x)]^*\otimes
F_x$ are mutually orthogonal in $C^*(W^u,F)$, and that if $x\in B$,
$f,\ f'\in F_x$,
\begin{align}\label{2.30}
\left\langle W^u( {x})^*\otimes f, W^u( {x})^*\otimes
f'\right\rangle =\left\langle f,f'\right\rangle_{b^{ {F}_{ {x}}}}.
\end{align}

We construct the equivariant symmetric bilinear form $b_{{\rm
det}(C^*(W^u,F),G)}$ on ${{\rm det}(C^*(W^u,F),G)}$ as in Definition
\ref{t2.1}.

\begin{defn}\label{t2.3}
The symmetric bilinear form on the determinant line of the
cohomology of the Thom-Smale cochain complex
$(C^*(W^u,F),\partial)$, in the sence of Definition \ref{t2.2}, is
called the equivariant Milnor symmetric bilinear torsion and is
denoted by $b^{{\cal M}, -X}_{{\rm det}(H^{*}(W^u,F),G)}$.
\end{defn}

Take $g\in G$. Set
\begin{align}\label{2.31}
M_g=\{x\in M,\ gx=x\}.
\end{align}
Since $G$ is a compact group, $M_g$ is a smooth compact submanifold
of $M$. Let $N$ be the normal bundle to $M_g$ in $M$.  By
\cite[Proposition 1.13]{BZ2}, we know that $f|_{M_g}$ is a Morse
function on $M_g$, and $X|_{M_g}$ is a smooth section of $TM_g$. For
$g\in G$, set
\begin{align}\label{2.32}
B_g=B\cap M_g.
\end{align}
Then $B_g$ is the set of critical points of $f|_{M_g}$.
\begin{defn}\label{t2.4}
If $x\in B_g$, let ${\rm ind}_g(x)$ be the index of $f|_{M_g}$ at
$x$.
\end{defn}

Let now $b^{F_x}$, $b'^{F_x}$ $(x\in B)$ be two $G$-invariant
nondegenerate symmetric bilinear forms on $F_x$. Let $b^{{\cal M},
-X}_{{\rm det}(H^{*}(W^u,F),G)}$, $b'^{{\cal M}, -X}_{{\rm
det}(H^{*}(W^u,F),G)}$ be the corresponding equivariant Milnor
symmetric bilinear torsions. By \cite[Theorem 1.15]{BZ2} and
\cite[Proposition 2.5]{SZ}, we have the following theorem.
\begin{thm}\label{t2.5}
For $g\in G$, the following identity holds
\begin{align}\label{2.33}
b'^{{\cal M}, -X}_{{\rm det}(H^{*}(W^u,F),G)}(g)=b^{{\cal M},
-X}_{{\rm det}(H^{*}(W^u,F),G)}(g)\prod_{x\in B_g}\exp\left({\rm
Tr}_{F_x}\left[g\log\left(b'^{F_x}\over
b^{F_x}\right)\right]\right)^{(-1)^{{\rm ind}_g(x)}}.
\end{align}
\end{thm}

\subsection{\normalsize  Equivariant Ray-Singer symmetric bilinear torsion }
\label{s2.4} We continue the discussion of the previous subsection.
However, we do not use the Morse function and make transversality
assumptions.

For any $0\leq i\leq n$, denote
\begin{align}\label{2.34}
\Omega^i( {M}, {F})=\Gamma\left(\Lambda^i(T^* {M})\otimes
 {F}\right),\ \ \ \
\Omega^*( {M}, {F})=\bigoplus_{i=0}^n\Omega^i( {M}, {F}).
\end{align}
Let $d^{ {F}}$ denote the natural exterior differential on
$\Omega^*( {M}, {F})$ induced from $\nabla^{ {F}}$ which maps each
$\Omega^i( {M}, {F})$, $0\leq i\leq n$, into $\Omega^{i+1}( {M},
{F})$.

The group $G$ acts naturally on $\Omega^*( {M}, {F})$. Namely, if
$g\in G$, $s\in \Omega^*( {M}, {F})$, set
$$gs(x)=g_* s(g^{-1}x),\ x\in M.$$

Let $g^F$ be a $G$-invariant Hermitian metric on $F$. The
$G$-invariant Riemannian metric $g^{TM}$ and $g^F$ determine a
natural inner product $\langle\ ,\ \rangle_g$ (that is, a
pre-Hilbert space structure) on $\Omega^*( {M}, {F})$  (cf.
\cite[(2.2)]{BZ1} and \cite[(2.3)]{BZ2}).

Let $d_g^{F_*}$ be the formal adjoint of $d^F$ with respect to
$\langle\ ,\ \rangle_g$ and $D_g=d^{F}+d_g^{F_*}$.

On the other hand  $g^{TM}$ and the $G$-invariant symmetric bilinear
form $b^F$ determine together a $G$-invariant symmetric bilinear
form on $\Omega^*( {M}, {F})$ such that if $u=\alpha f$, $v=\beta
g\in\Omega^*( {M}, {F})$ such that $\alpha,\ \beta\in \Omega^*( {M}
)$, $f,\ g\in \Gamma(F)$, then
\begin{align}\label{2.35}
 \langle u,v\rangle_b=\int_M (\alpha\wedge * \beta) b^F(f,g),
\end{align}
where $*$ is the Hodge star operator (cf. \cite{Z}).

Consider the   de Rham complex
\begin{multline}\label{2.36}
\left( \Omega^*( {M}, {F}) ,d^{ {F}}\right):0\rightarrow
 \Omega^0( {M}, {F}) \stackrel{d^{ {F}}}{\rightarrow}
 \Omega^1( {M}, {F})\rightarrow \cdots\\
\stackrel{d^{ {F}}}{\rightarrow}
 \Omega^n( {M}, {F} )\rightarrow 0.
\end{multline}

Let $d^{F *}_b: \Omega^*( {M}, {F})\rightarrow \Omega^*( {M}, {F})$
denote the formal adjoint of $d^{ {F}}$ with respect to
$G$-invariant the symmetric bilinear form in (\ref{2.35}). That is,
for any $u,\ v\in\Omega^*( {M}, {F})$, one has
\begin{align}\label{2.37}
 \left\langle d^Fu,v\right\rangle_b= \left\langle
 u,d^{F*}_bv\right\rangle_b.
\end{align} Set
\begin{align}\label{2.38}
 {D}_b=d^{ {F}}+d^{ {F}*}_b,\ \ \  {D}_b^2=\left(d^{
{F}}+d^{ {F}*}_b\right)^2=d^{ {F}*}_bd^{ {F}}+ d^{ {F}}d^{ {F}*}_b.
\end{align}
Then the Laplacian $ {D}^2_b$ preserves the ${\bf Z}$-grading of
$\Omega^*( {M}, {F})$.

As was pointed out in \cite{BH1} and \cite{BH2}, $D_b^2$ has the
same principal symbol as the usual Hodge Laplacian (constructed
using the inner product on $\Omega^*(M,F)$ induced from
$(g^{TM},g^F)$) studied for example in \cite{BZ1}.

We collect some well-known facts concerning $D_b^2$ as in
\cite[Proposition 4.1]{BH2}, where the reference \cite{S} is
indicated.

\begin{prop}\label{t2.6} The following properties hold for the Laplacian
$D_b^2$:

(i) The spectrum of $D_b^2$ is discrete. For every $\theta>0$ all
but finitely many points of the spectrum are contained in the angle
$\{ z\in{\bf C}|-\theta<{\rm arg}(z)<\theta\}$;

(ii) If $\lambda$ is in the spectrum of $D_b^2$, then the image of
the associated spectral projection is finite dimensional and
contains  smooth forms only. We refer to this image as the
(generalized) $\lambda$-eigen space of $D_b^2$ and denote it  by
$\Omega^*_{\{\lambda\}}(M,F)$. There exists $N_\lambda\in{\bf N}$
such that
\begin{align}\label{2.39}
  \left.\left(D_b^2-\lambda\right)^{N_\lambda}\right|_{\Omega^*_{\{\lambda\}}(M,F)}=0.
\end{align} We have a $D_b^2$-invariant  $\langle\ ,\
\rangle_b$-orthogonal decomposition
\begin{align}\label{2.40}
  \Omega^* (M,F)=\Omega^*_{\{\lambda\}}(M,F)\oplus \Omega^*_{\{\lambda\}}(M,F)^{\perp } .
\end{align} The restriction of $D_b^2-\lambda$ to $\Omega^*_{\{\lambda\}}(M,F)^{\perp }$ is invertible;

(iii) The decomposition (\ref{2.40}) is invariant under $d^F$ and
$d^{F*}_b$;

(iv) For $\lambda\neq \mu$, the eigen spaces
$\Omega^*_{\{\lambda\}}(M,F)$ and $\Omega^*_{\{\mu\}}(M,F)$ are
$\langle\ ,\ \rangle_b$-orthogonal to each other.
\end{prop}

For any $a\geq 0$, set
\begin{align}\label{2.41}
  \Omega^*_{[0,a]} (M,F)=\bigoplus_{0\leq |\lambda|\leq a}\Omega^*_{\{\lambda\}}(M,F).
\end{align} Let $\Omega^*_{[0,a]} (M,F)^\perp$ denote the $\langle\ ,\ \rangle_b$-orthogonal
complement to $\Omega^*_{[0,a]} (M,F)$. Obviously, each
$\Omega^*_{\{\lambda\}}(M,F)$ is a $G$-invariant subspace.

By \cite[(29)]{BH2} and Proposition \ref{t2.6}, one sees that
$(\Omega^*_{[0,a]} (M,F),d^F)$ forms a finite dimensional complex
whose cohomology equals to that of $(\Omega^*  (M,F),d^F)$.
Moreover, the $G$-invariant symmetric bilinear form $\langle\ ,\
\rangle_b$ clearly induces a nondegenerate $G$-invariant symmetric
bilinear form on each $\Omega^i_{[0,a]} (M,F)$ with $0\leq i\leq n$.
By Definition \ref{t2.2} one then gets a symmetric bilinear torsion
$b^{\rm RS}_{\det (H^*(\Omega^*_{[0,a]} (M,F)),G)}$ on $\det
H^*(\Omega^*_{[0,a]} (M,F),d^F)=\det H^*(\Omega^*  (M,F),d^F)$.

For any $0\leq i\leq n$, let $D_{b,i}^2$ be the restriction of
$D_{b}^2$ on $\Omega^i(M,F)$. Then it is shown in \cite{BH2} (cf.
\cite[Theorem 13.1]{S}) that for any $a\geq 0$, $g\in G$ the
following is well-defined,
\begin{align}\label{2.42}
   {\det}'\left(D^2_{b,(a,+\infty),i}\right)(g)
   =\exp\left(-\left.{\partial \over\partial s}\right|_{s=0}{\rm
   Tr}\left[g
   \left(\left. D^2_{b,i}\right|_{\Omega^*_{[0,a]} (M,F)^\perp}\right)^{-s}\right]
   \right).
\end{align}
\begin{defn}\label{t2.7}
If $g\in G$, set
\begin{align}\label{2.43}
b^{\rm RS}_{{\rm det}(H^*(M,F),G)}(g)=b^{\rm RS}_{\det
(H^*(\Omega^*_{[0,a]} (M,F)),G)}(g)\prod_{i=0}^n \left(
{\det}'\left(D^2_{b,(a,+\infty),i}\right)(g)\right)^{(-1)^ii},
\end{align}
by \cite[Proposition 4.7]{BH2}, we know that $b^{\rm RS}_{{\rm
det}(H^*(M,F),G)}$ does not depend on the choice of $a\geq 0$, and
is called the equivariant Ray-Singer symmetric bilinear torsion on
$\det H^*(\Omega^* (M,F),d^F)$.
\end{defn}

\subsection{\normalsize An anomaly formula for the equivariant Ray-Singer symmetric bilinear torsion }
\label{s2.5}

We continue the discussion of the above subsection.
\begin{defn}\label{t2.8}
Let $\theta_g(F,b^F)$ be the $1$-form on $M_g$
\begin{align}\label{2.44}
\theta_g(F,b^F)={\rm Tr}\left[g(b^F)^{-1}\nabla^F b^F\right].
\end{align}
\end{defn}

Clearly $M_g$ is a totally geodesic submanifold of $M$. Let
$g^{TM_g}$ be the Riemannian metric induced by $g^{TM}$ on $TM_g$.
Let $\nabla^{TM_g}$ be the Levi-Civita connection on $(TM_g,
g^{TM_g})$.

Let $e(TM_g,\nabla^{TM_g})$ be the Chern-Weil representative of the
rational Euler class of $TM_g$, associated to the metric preserving
connection $\nabla^{TM_g}$. Then
\begin{align}\label{2.45}
e(TM_g,\nabla^{TM_g})={\rm Pf}\left[{R^{TM_g}}\over{2\pi}\right]
{\rm if}\ {\rm dim}M_g\ {\rm is\ even},
\end{align}
$$\ \ \ \ \ \ \ \ \ \ \ \ \ \ \ \ \ \ \ \ \ \ \ \ \ \ \ \ \ \ \ \ \ 0\ \ \ \ \ {\rm if}\ {\rm dim}M_g\ {\rm is\ odd}.$$

Let $g'^{TM}$ be another $G$-invariant metric and let
$\nabla'^{TM_g}$ be the corresponding Levi-Civita connection on
$TM_g$. Let $\widetilde{e}(TM_g,\nabla^{TM_g},\nabla'^{TM_g})$ be
the Chern-Simons class of ${\rm dim}M_g-1$ forms on $M_g$, such that
\begin{align}\label{2.46}
d\widetilde{e}\left(TM_g,\nabla^{TM_g},\nabla'^{TM_g}\right)=e\left(TM_g,\nabla'^{TM_g}\right)
-e\left(TM_g,\nabla^{TM_g}\right).
\end{align}

Let $b'^F$ be another $G$-invariant nondegenerate symmetric bilinear
form on $F$.

Let $b'^{\rm RS}_{{\rm det}(H^*(M,F),G)}$ denote the equivariant
Ray-Singer symmetric bilinear torsion associated to $g'^{TM}$ and
$b'^F$.

By \cite[Remark 6.4]{SZ} and \cite[Theorem 2.7]{BZ2}, we have the
following extension of the anomaly formula of \cite[Theorem
2.9]{SZ}.
\begin{thm}\label{t2.9}
If $b^F$, $b'^F$ lie in the same homotopy class of nondegenerate
symmetric bilinear forms on $F$, then for $g\in G$ the following
identity holds,
\begin{multline}\label{2.47}
\left({b'^{\rm RS}_{{\rm det}(H^*(M,F),G)}\over{b^{\rm RS}_{{\rm
det}(H^*(M,F),G)}}}\right)(g)=\exp\left(\int_{M_g}{\rm
Tr}\left[g\log\left({b'^F}\over{b^F}\right)\right]e\left(TM_g,\nabla^{TM_g}\right)\right)\\
\cdot\exp\left(-\int_{M_g}\theta_g
(F,b'^F)\widetilde{e}\left(TM_g,\nabla^{TM_g},\nabla'^{TM_g}\right)\right).
\end{multline}
\end{thm}
{\it Proof.} Let $b_{l}^{F}$ is a smooth one-parameter family of
fiber wise non-degenerate symmetric bilinear forms on $F$ and
$(g_{l}^{TM}, g_{l}^{F})$ be a smooth family of metrics on $TM, F$.

By \cite[(6.4)]{SZ}, we have
\begin{multline}\label{2.1*}
e^{-tD_b^2}=e^{-tD_g^2}+\sum_{k=1}^{n}(-1)^kt^k\int_{\Delta_k}e^{-t_1tD_g^2}B_{b,g}
e^{-t_2tD_g^2}\cdots B_{b,g}e^{-t_{k+1}tD_g^2}dt_1\cdots dt_k\\
+ (-1)^{n+1}t^{n+1}\int_{\Delta_{n+1}}e^{-t_1tD_g^2}B_{b,g}
e^{-t_2tD_g^2}\cdots B_{b,g}e^{-t_{n+2}tD_b^2}dt_1\cdots dt_{n+1},
\end{multline}
where $\Delta_k$, $1\leq k\leq n+1$, is the $k$-simplex defined by
$t_1+\cdots+t_{k+1}=1$, $t_1\geq 0$, $\cdots, $ $t_{k+1}\geq 0$ and
$B_{b,g}$ is defined in \cite[(6.3)]{SZ}.

Proceeding as in \cite[Section 4]{BZ1}, we first calculate the
asymptotics as $t\to 0$ of ${\rm Tr}_{s}[g(b_{l}^{F})^{-1}{{\partial
b_{l}^{F}}\over{\partial l}}\exp(-tD_{b_l}^2)]$. Here the metric
$g^{TM}$ will be fixed.

By the same proof in \cite[proposition 6.1]{SZ}, we have that as
$t\to 0^{+}$,
\begin{align}\label{2.2*}
t^{n+1}\int_{\Delta_{n+1}}{\rm Tr}_s\left[g
(b_{l}^{F})^{-1}{{\partial b_{l}^{F}}\over{\partial l}}
e^{-t_1tD_g^2}B_{b,g} e^{-t_2tD_g^2}\cdots
B_{b,g}e^{-t_{n+2}tD_b^2}\right]dt_1\cdots dt_{n+1}\rightarrow 0.
\end{align}
Also, by \cite[(6.22)]{SZ}, we have that for $ 1\leq k\leq n$,
$(t_1,\cdots,t_{k+1})\in\Delta_k$,
\begin{align}\label{2.3*}
\lim_{t\rightarrow 0^+}t^k {\rm Tr}_s\left[g
(b_{l}^{F})^{-1}{{\partial b_{l}^{F}}\over{\partial l}}
e^{-t_1tD_g^2}B_{b,g} e^{-t_2tD_g^2}\cdots
B_{b,g}e^{-t_{k+1}tD_g^2}\right]=0.
\end{align}
So that by (\ref{2.1*})-(\ref{2.3*}) we have that
\begin{align}\label{2.4*}
\lim_{t\to 0}{\rm Tr}_{s}\left[g(b_{l}^{F})^{-1}{{\partial
b_{l}^{F}}\over{\partial
l}}\exp\left(-tD_{b_l}^2\right)\right]=\lim_{t\to 0}{\rm
Tr}_{s}\left[g(b_{l}^{F})^{-1}{{\partial b_{l}^{F}}\over{\partial
l}}\exp\left(-tD_{g_l}^2\right)\right].
\end{align}

Now we assume that the nondegenerate symmetric bilinear form on $F$
is fixed, and the metric $g_{l}^{TM}$ on $TM$ depends on $l$.

Let $*_{l}$ be the Hodge star operator associated to $g_{l}^{TM}$.

By \cite[(4.70), (4.74)]{BZ1}, analogues of \cite[(6.5), (6.24),
(6.26), (6.27)]{SZ} replacing $N$ by
$*_{l}^{-1}{{\partial*_{l}^{TM}}\over{\partial l}}$ and
\cite[Chapter 6]{BGV}, we have that
\begin{multline}\label{2.5*}
\lim_{t\to 0}{\rm Tr}_{s}\left[g*_{l}^{-1}{{\partial*_{l}^{TM}}\over
{\partial l}}\exp\left(-tD_{b_l}^2\right)\right]=\lim_{t\to 0}{\rm
Tr}_{s}\left[g*_{l}^{-1}{{\partial*_{l}^{TM}}\over {\partial
l}}\exp\left(-tD_{g_l}^2\right)\right]\\
-{1\over 2}\int_{M_g} \int^B {\rm Tr}\left[g\left(  \sum_{i,\,
j=1}^n
   e_i \wedge\widehat{e_j}\left(\nabla_{e_i}^u\omega^F\left(e_j\right)\right)+{1\over
  2}\left[
   \omega^F ,\widehat{ \omega^F_g}-\widehat{ \omega^F }
   \right]\right)\right]\\
   \cdot\left(-\sum_{1\leq i,j\leq n}{1\over 2}\left\langle\left(g_{l}^{TM}\right)^{-1}{{\partial g_{l}^{TM}}\over {\partial l}}e_{i},e_{j}\right\rangle_{g_{l}^{TM}}e_{i}\wedge\widehat{e_{j}}\right)\exp\left(-{\dot{R_l}^{TM_g}\over
   2}\right)\\
   =\lim_{t\to 0}{\rm
Tr}_{s}\left[g*_{l}^{-1}{{\partial*_{l}^{TM}}\over {\partial
l}}\exp\left(-tD_{g_l}^2\right)\right]-{1\over
2}\int_{M_g}\int^{B}\nabla_l^{TM}\varphi{\rm
Tr}\left[g\omega^{F}\right]\\
\cdot\left(-\sum_{1\leq i,j\leq n}{1\over
2}\left\langle\left(g_{l}^{TM}\right)^{-1}{{\partial
g_{l}^{TM}}\over {\partial
l}}e_{i},e_{j}\right\rangle_{g_{l}^{TM}}e_{i}\wedge\widehat{e_{j}}\right)\exp\left(-{\dot{R_l}^{TM_g}\over
   2}\right)\\
=\int_{M_g}\left\{\int^{B}\nabla_{l}^{TM}\left({1\over 4}\sum_{1\leq
i,j\leq n}\left\langle\left(g_{l}^{TM}\right)^{-1}{{\partial
g_{l}^{TM}}\over {\partial
l}}e_{i},e_{j}\right\rangle_{g_{l}^{TM}}e_{i}\wedge\widehat{e_{j}}\right)\right.\\
\left.\cdot\exp\left(-{\dot{R_l}^{TM_g}\over
   2}\right)\wedge\varphi\theta_{g}\left(F,b^F\right)\right\}.
\end{multline}

From (\ref{2.4*}), (\ref{2.5*}) and the calculations in
\cite[Section 4]{BZ1}, we get (\ref{2.47}).

The proof of Theorem \ref{t2.9} is completed. \ \ Q.E.D.

\section{A formula relating equivariant Milnor and equivariant Ray-Singer symmetric bilinear torsions} \label{s3}
\setcounter{equation}{0}

In this section, we state the main result of this paper, which is an
explicit comparison result between the equivariant Milnor symmetric
bilinear torsion and equivariant Ray-Singer symmetric bilinear
torsion.

We assume that we are in the same situation as in Sections
 \ref{s2.2}-\ref{s2.4}.
By a simple argument of Helffer-Sj\"ostrand \cite[Proposition
5.1]{HS} (cf. \cite[Section 7b)]{BZ1}), we  may and we well assume
that $g^{TM}$ there satisfies the following property  without
altering the  Thom-Smale cochain complex $(C^*(W^u, {F}),
{\partial})$,

 (*): For any $x\in B$, there is a system of coordinates
$y=(y^1,\cdots,y^n)$ centered at $x$ such that near $x$,
\begin{align}\label{2.48}
g^{TM}=\sum_{i=1}^n \left|dy^i\right|^2,\ \ \ f(y)=f(x)-{1\over
2}\sum_{i=1}^{{\rm ind}(x)}\left|y^i\right|^2+{1\over 2}\sum_{i={\rm
ind}(x)+1}^{n}\left|y^i\right|^2.
\end{align}

By a result of Laudenbach \cite{La}, $\{W^u(x):x\in B\}$ form a CW
decomposition of $M$.

For any $ {x}\in  {B}$, $ {F}$ is canonically trivialized over each
cell $W^u( {x})$.

Let $ {P}_\infty$ be the de Rham map defined by
\begin{align}\label{2.49}
\alpha\in\Omega^*( {M}, {F})
 \rightarrow
 {P}_\infty\alpha=\sum_{ {x}\in
 {B}}W^u( {x})^*\int_{W^u( {x})}\alpha\in
C^*(W^u, {F}).
\end{align}

 By the Stokes theorem,   one has
\begin{align}\label{2.50}
 {\partial} {P}_\infty= {P}_\infty
d^{ {F}}.
\end{align}
Moreover, it is shown in \cite{La} that $P_\infty$ is a  ${\bf
Z}$-graded   quasi-isomorphism, inducing a canonical  isomorphism
\begin{align}\label{2.51}
 {P}^{ H}_\infty:{  H}^*
 \left(\Omega^*( {M}, {F}),d^{ {F}}\right)\rightarrow
{  H}^*\left(C^*\left(W^u, {F}\right), {\partial}\right),
\end{align}
which in turn induces a natural isomorphism between the determinant
lines,
\begin{align}\label{2.52}
 {P}^{\det {  H}}_\infty:\det{  H}^*\left(\Omega^*\left( {M}, {F}\right),d^{ {F}}
 \right)\rightarrow
\det{ H}^*\left(C^*\left(W^u, {F}\right), {\partial}\right).
\end{align}
Also by \cite[Theorem 1.11]{BZ2}, we know that ${P}_\infty$ commutes
with $G$, and $P^{H}_\infty$ is the canonical identification of the
corresponding cohomology groups as $G$-spaces.

Now let $h^{TM}$ be an arbitrary smooth metric on $TM$.

By Definition \ref{t2.7}, one has an associated  equivariant
Ray-Singer symmetric bilinear torsion $ b^{\rm RS}_{{\rm
det}(H^*(M,F),G )} $ on $\det{ H}^* (\Omega^* ( {M}, {F} ),d^{ {F}}
  )$.
From (\ref{2.52}), one gets a well-defined equivariant symmetric
bilinear form
\begin{align}\label{2.53}
 {P}^{\det {  H}}_\infty \left(b^{\rm RS}_{{\rm
det}(H^*(M,F),G )} \right)
\end{align}
on $\det{ H}^* (C^* (W^u, {F} ), {\partial} )$.

On the other hand, by Definition \ref{t2.3}, one has a well-defined
equivariant Milnor symmetric bilinear torsion
 $
b^{{\cal M},-X}_{{\rm det}(H^*(M,F),G)}$ on $\det{ H}^*(C^*(W^u,
{F}), {\partial})$,
 where $X=\nabla f$ is the gradient vector field of $f$ associated
to $g^{TM}$.

Let $M_g=\cup_{j=1}^{m}M_{g,j}$ be the decomposition of $M_g$ into
its connected components. Clearly ${\rm Tr}_{F}[g]$ is constant on
each $M_{g,j}$.

Let $N$ be the normal bundle to $M_g$ in $M$. We identify $N$ to the
orthogonal bundle to $TM_g$ in $TM|_{M_g}$.

Take $x\in B_g$. Then $g$ acts on $T_x M$ as a linear isometry. Also
$$TM_g=\{Y\in TM|_{M_g},\ gY=Y\}.$$
Moreover $g$ acts on $N$. Let $e^{\pm i\beta_1},\cdots,e^{\pm
i\beta_{q}}(0<\beta_j\leq\pi)$ be the locally constant distinct
eigenvalues of $g|_N$. Then $N$ splits orthogonally as
$$N=\bigoplus_{j=1}^{q}N^{\beta_j}.$$
For $1\leq j\leq q$, $g$ acts on $N^{\beta_j}$ as an isometry, with
eigenvalues $e^{\pm i\beta_j}$. In particular, if $e^{\pm
i\beta_j}\neq -1$, $N^{\beta_j}$ is even dimensional.

Take $x\in B_g$. Since $f$ is $g$-invariant, $d^{2}f(x)$ is also
$g$-invariant. Therefore the decomposition
$$T_x M=T_x M_g\oplus\bigoplus_{j=1}^{q}N^{\beta_j}$$
is orthogonal with respect to $d^2 f(x)$. On $T_x M_g$, the index of
$d^2 f(x)|_{T_x M_g\times T_x M_g}$ was already denoted ${\rm ind}_g
(x)$. Let $n_{+}(\beta_j)(x)$ (${\rm resp.}\ n_{-}(\beta_j)(x)$) be
the number of positive (resp. negative) eigenvalues of $d^2
f(x)|_{N^{\beta_j}}$. Then if $e^{\pm i\beta_j}\neq -1,
n_{\pm}(\beta_j)(x)$ is even.

Let $\psi(TM_g,\nabla^{TM_g})$ be the Mathai-Quillen current
(\cite{MQ}) over $TM_g$, associated to $h^{TM}$, defined  in
\cite[Definition 3.6]{BZ1}. As indicated in \cite[Remark 3.8]{BZ1},
the pull-back current $X^*\psi(TM_g,\nabla^{TM_g})$ is well-defined
over $M_g$.

The main result of this paper, which generalizes \cite[Theorem
3.1]{SZ} to the equivariant case.

\begin{thm}\label{t2.10}
For $g\in G$, the following identity in $\bf C$ holds,
\begin{multline}\label{2.54}
{{P}^{\det {  H}}_\infty \left(b^{\rm RS}_{{\rm det}(H^*(M,F),G )}
\right)\over {b_{{\rm det}(H^*(W^u,F),G)}^{{\cal
M},-X}}}(g)=\exp\left(-\int_{M_g}\theta_g(F,b^F)X^*\psi\left(TM_g,\nabla^{TM_g}\right)\right)\\
\cdot\exp\left(-{1\over 4}\sum_{x\in B_g}(-1)^{{\rm ind}_g
(x)}\sum_{j}(n_{+}(\beta_j)(x)-n_{-}(\beta_j)(x))\right.\\
\left.\cdot\left({\Gamma'\over\Gamma}\left({\beta_j\over
2\pi}\right)+{\Gamma'\over\Gamma}\left(1-{\beta_j\over
2\pi}\right)-2\Gamma'(1)\right)\cdot{\rm
Tr}\left[g|_{F_x}\right]\right).
\end{multline}
\end{thm}

\begin{Rem}\label{t2.11} By proceeding similarly as in
\cite[Section 5b)]{BZ2}, in order to prove (\ref{2.54}), we may well
assume that $h^{TM}=g^{TM}$. Moreover, we may assume that
 $b^F$, as well as the Hermitian metric $h^F$ on $F$, are
  flat on an open neighborhood of the zero set $B$ of
$X$. From now on, we will make these assumptions.
\end{Rem}

\section{A proof of Theorem 3.1} \label{s3}
\setcounter{equation}{0}

We assume that the assumptions in Remark \ref{t2.11} hold.

For any $T\in {\bf R}$, let $b^F_T$ be the deformed symmetric
bilinear form on $F$ defined by
\begin{align}\label{3.1}
  b^F_T(u,v)=e^{-2Tf}b^F(u,v).
\end{align}
Let $d^{F*}_{b_T}$ be the associated formal adjoint in the sense of
(\ref{2.37}). Set
\begin{align}\label{3.2}
 {D}_{b_T}=d^{ {F}}+d^{ {F}*}_{b_T},\ \ \  {D}_{b_T}^2=\left(d^{
{F}}+d^{ {F}*}_{b_T}\right)^2=d^{ {F}*}_{b_T}d^{ {F}}+ d^{ {F}}d^{
{F}*}_{b_T}.
\end{align}

Let $\Omega^*_{[0,1],T}(M,F)$ be defined as in (\ref{2.41}) with
respect to $ {D}_{b_T}^2$, and let $\Omega^*_{[0,1],T}(M,F)^\perp$
be the corresponding $\langle\ ,\ \rangle_{b_T}$-orthogonal
complement.

Let $P_T^{[0,1]}$ be the orthogonal projection from $\Omega^*(M,F)$
to $\Omega^*_{[0,1],T}(M,F)$ with respect to the inner product
determined by $g^{TM}$ and $g_T^F=e^{-2Tf}g^F$.  Set
$P_T^{(1,+\infty)}={\rm Id}-P_T^{[0,1]}$.

Following \cite[(5.9)-(5.10)]{BZ2}, we introduce the notations
\begin{align}\label{3.3}
\chi_{g}(F)=\sum_{j}{\rm Tr}_{F|_{M_{g,j}}}[g]\sum_{x\in B_g\cap
M_{g,j}}(-1)^{{\rm ind}_{g}(x)},
\end{align}
$$\widetilde{\chi}'_g(F)=\sum_{j}{\rm Tr}_{F|_{M_{g,j}}}[g]\sum_{x\in B_g\cap
M_{g,j}}(-1)^{{\rm ind}_{g}(x)}{\rm ind}(x),$$
$${\rm Tr}_{s}^{B_g}[f]=\sum_{j}{\rm Tr}_{F|_{M_{g,j}}}[g]\sum_{x\in B_g\cap
M_{g,j}}(-1)^{{\rm ind}_{g}(x)}f(x).$$

Let $N$ be the number operator on $\Omega^*(M,F)$ acting on
$\Omega^i(M,F)$ by multiplication by $i$.

By the technique developed in \cite{SZ} and the corresponding
results in \cite{BZ2}, we easily get the following intermediate
results. The sketch of the proofs will be outlined in Section 5.

\begin{thm}\label{t3.1} {\rm (Compare with \cite[Theorem 5.5]{BZ2} and \cite[Theorem 3.3]{SZ})}
 Let $P_T^{[0,1]}$ be the restriction of $P_\infty$ on
$\Omega^*_{[0,1],T}(M,F)$, let $P_T^{[0,1],\det H}$ be the induced
isomorphism on cohomology, then the following identity holds,
\begin{align}\label{3.4}
\lim_{T\rightarrow +\infty} {P_T^{[0,1],\det H}\left(b^{\rm
RS}_{\det (H^*(\Omega^*_{[0,1],T} (M,F)),G)}\right)\over b_{{\rm
det}(H^*(W^u,F),G)}^{{\cal M},-X} }(g)\left({T\over
\pi}\right)^{{n\over 2}\chi_g(F)-
\widetilde{\chi}'_g(F)}\exp\left(2{\rm Tr}^{B_g}_s[f]T\right)
\end{align}$$ =1.$$
\end{thm}

\begin{thm}\label{t3.2} {\rm (Compare with \cite[Theorem 5.7]{BZ2} and \cite[Theorem 3.4]{SZ})} For any $t>0$,
\begin{align}\label{3.5}
\lim_{T\rightarrow +\infty}  {\rm Tr}_s\left[g
N\exp\left(-tD_{b_T}^2\right)P_T^{(1,+\infty)}\right]=0.
\end{align}
Moreover, for any  $d>0$ there exist $c>0$, $C>0$ and $T_0\geq 1$
such that for any $t\geq d$ and  $T\geq T_0$,
\begin{align}\label{3.6}
 \left|{\rm Tr}_s\left[gN\exp\left(-tD_{b_T}^2\right)P_T^{(1,+\infty)}\right]\right|\leq c\exp(-Ct).
\end{align}
\end{thm}

\begin{thm}\label{t3.3} {\rm (Compare with \cite[Theorem 5.8]{BZ2} and \cite[Theorem 3.5]{SZ})}
For $T\geq 0$ large enough, then
\begin{align}\label{3.7}
  \lim_{T\to +\infty}{\rm
  Tr}_s\left[gNP_{T}^{[0,1]}\right]=\widetilde{\chi}'_{g}(F).
\end{align}
Also,
\begin{align}\label{3.8}
 \lim_{T\rightarrow +\infty}  {\rm Tr}\left[D_{b_T}^2P_T^{[0,1]}\right]=0 .
\end{align}
\end{thm}

For the next results, we will make use the same notation for
Clifford multiplications and Berezin integrals as in \cite[Section
4]{BZ1}.

\begin{thm}\label{t3.4} {\rm (Compare with \cite[Theorem 5.9]{BZ2} and \cite[Theorem 3.6]{SZ})}
As $t\rightarrow 0$, the following identity holds,
\begin{align}\label{3.9}
   {\rm Tr}_s\left[gN\exp\left(-tD_{b_T}^2\right)\right]={n\over 2}\chi_g(F)+O(t)\ \ if\ n\ {\rm is\ even},
\end{align}
$$
  =\int_{M_g}{\rm Tr}_F[g]\int^BL\exp\left(-{\dot{R}^{TM_g}\over 2}\right)
  {1\over \sqrt{t}}+O\left(\sqrt{t}\right)\ \ if\ n \ {\rm is\ odd}.
$$
\end{thm}

\begin{thm}\label{t3.5} {\rm (Compare with \cite[Theorem A.1]{BZ2} and \cite[Theorem 3.7]{SZ})}
There exist $0<\alpha\leq 1$, $C>0$ such that for any $0<t\leq
\alpha$, $0\leq T\leq {1\over t}$, then
\begin{multline}\label{3.10}
 \left|  {\rm Tr}_s\left[gN\exp\left(-\left(tD_{b }+T\widehat{c}
 (\nabla f)\right)^2\right)\right]-
 {1\over t}\int_{M_g}{\rm Tr}_F[g]\int^BL\exp\left(-B_{T^2}\right)\right.\\
\left.  -{T\over
2}\int_{M_g}\theta_g\left(F,b^F\right)\int^B\widehat{df}\exp\left(-B_{T^2}\right)-
  {n\over 2}\chi_g(F)\right|\leq Ct.
\end{multline}
\end{thm}

\begin{thm}\label{t3.6} {\rm (Compare with \cite[Theorem A.2]{BZ2} and \cite[Theorem 3.8]{SZ})} For any
$T>0$, the following identity holds,
\begin{multline}\label{3.11}
 \lim_{t\rightarrow 0}  {\rm Tr}_s\left[gN\exp\left(-\left(tD_{b }+{T\over t}\widehat{c}
 (\nabla f)\right)^2\right)\right] =
 \sum_{j}{\rm Tr}\left[g|_{F|_{M_{g,j}}}\right]\\
 \cdot\left({1\over{1-e^{-2T}}}
 \left(\left(1+e^{-2T}\right)\sum_{x\in{B\cap{M_{g,j}}}}(-1)^{{\rm ind}_g(x)}{{\rm ind}_g}(x)
 -{\rm dim}M_{g,j}e^{-2T}\chi(M_{g,j})\right)\right)\\
 +\sum_{j}{\rm Tr}\left[g|_{F|_{M_{g,j}}}\right]\sum_{k}{\sinh(2T)\over{\cosh(2T)-\cos(\beta_k)}}\sum_{x\in
 {B\cap{M_{g,j}}}}(-1)^{{\rm ind}_g(x)}n_{-}({\beta_k})(x)\\
 -\sum_{j}{\rm Tr}\left[g|_{F|_{M_{g,j}}}\right]\sum_{k}{1\over
 2}\left({\sinh(2T)\over{\cosh(2T)-\cos(\beta_k)}}-1\right){\rm
 dim}N^{\beta_k}\chi(M_{g,j}).
\end{multline}
\end{thm}

\begin{thm}\label{t3.7} {\rm (Compare with \cite[Theorem A.3]{BZ2} and \cite[Theorem 3.9]{SZ})}
There exist $\alpha\in(0,1]$, $c>0$, $C>0$ such that for any $t\in
(0,\alpha]$, $T\geq 1$, then
\begin{align}\label{3.12}
 \left|  {\rm Tr}_s\left[gN\exp\left(-\left(tD_{b }+{T\over t}\widehat{c}
 (\nabla f)\right)^2\right)\right]-
  \widetilde{\chi}'_g(F)\right|\leq c\exp(-CT)
  .
\end{align}
\end{thm}

Clearly, we may and we will assume that the number $\alpha>0$ in
Theorems \ref{t3.5} and \ref{t3.7} have been chosen to be the same.

Next, we use above theorems to give a proof of Theorem \ref{t2.10}.
Since the process is similar to it in \cite{SZ}, so we refer to it
for more details.

First of all, by the anomaly formula (\ref{2.47}), for any $T\geq
0$, $g\in G$, one has

\begin{multline}\label{4.13}
{P_T^{[0,1],\det H}\left(b^{\rm RS}_{\det
H^*(\Omega^*_{[0,1],T}(M,F),G)}\right)\over b_{{\rm
det}\left(H^*(W^u,F),G\right)}^{{\cal
M},-X}}(g)\\
\cdot\prod_{i=0}^n\left(\det\left(\left.
D^2_{b_T}\right|_{\Omega^*_{[0,1],T}(M,F)^\perp\cap
\Omega^{i}(M,F)}\right)(g)\right)^{(-1)^ii}\\
={P_\infty^{\det H}\left(b^{\rm RS}_{{\rm
det}(H^*(M,F),G)}\right)\over b_{{\rm det}(H^*(W^u,F),G)}^{{\cal
M},-X}}(g) \exp\left(-2T\int_{M_g}{\rm
Tr}_{F}[g]fe\left(TM_g,\nabla^{TM_g}\right)\right).
\end{multline}

From now on, we will write $a\simeq b$ for $a,\ b\in{\bf C}$ if
$e^a=e^b$. Thus, we can rewrite (\ref{4.13}) as

\begin{multline}\label{4.14}
\log\left({P_\infty^{\det H}\left(b^{\rm RS}_{{\rm
det}(H^*(M,F),G)}\right)\over b_{{\rm det}(H^*(W^u,F),G)}^{{\cal
M},-X}}(g) \right)\simeq\log\left({P_T^{[0,1],\det H}\left(b^{\rm
RS}_{\det H^*(\Omega^*_{[0,1],T}(M,F),G)}\right)\over b_{{\rm
det}\left(H^*(W^u,F),G\right)}^{{\cal M},-X}}(g)\right)\\
+\sum_{i=0}^{n}(-1)^{i} i\log\left(\det\left(\left.
D^2_{b_T}\right|_{\Omega^*_{[0,1],T}(M,F)^\perp\cap
\Omega^{i}(M,F)}\right)(g)\right)\\
+2T\int_{M_g}{\rm Tr}_{F}[g]fe\left(TM_g,\nabla^{TM_g}\right).
\end{multline}

Let $T_0>0$ be as in Theorem \ref{t3.2}. For any $T\geq T_0$
 and $s\in{\bf C}$ with ${\rm Re}(s)\geq n+1$, set
 \begin{align}\label{4.15}
\theta_{g,T}(s)={1\over \Gamma(s)}\int_0^{+\infty}t^{s-1}{\rm
 Tr}_s\left[gN\exp\left(-tD_{b_T}^2\right)P^{(1,+\infty)}_T\right]dt.
 \end{align}
By (\ref{3.6}), $\theta_{g,T}(s)$ is well defined and can be
extended to a meromorphic function which is holomorphic at $s=0$.
Moreover,
\begin{align}\label{4.16}
\sum_{i=0}^{n}(-1)^{i} i\log\left(\det\left(\left.
D^2_{b_T}\right|_{\Omega^*_{[0,1],T}(M,F)^\perp\cap
\Omega^{i}(M,F)}\right)(g)\right)\simeq-\left.{{\partial\theta_{g,T}(s)}\over{\partial
s}}\right|_{s=0}.
\end{align}

Let $d=\alpha^2$ with $\alpha$ being as in Theorem \ref{t3.7}. From
(\ref{4.15}) and Theorems \ref{t3.2}-\ref{t3.4}, one finds that
\begin{multline}\label{4.17}
\lim_{T\to +\infty}\left.{\partial\theta_{g,T}(s)\over{\partial
s}}\right|_{s=0}=\lim_{T\to +\infty}\int_{0}^{d}\left(
 {\rm
 Tr}_s\left[gN\exp\left(-tD_{b_T}^2\right)\right]-
 {a_{-1}\over\sqrt{t}}-{n\over 2}\chi_g(F) \right) {dt\over t}\\
- {2a_{-1}\over{\sqrt{d}}} -\left(\Gamma'(1)-{\log
d}\right)\left({n\over 2}\chi_g(F)-  \widetilde{\chi}'_g(F)\right).
\end{multline}

To study the first term in the right hand side of (\ref{4.17}), we
observe first that for any $T\geq 0$, one has
\begin{align}\label{4.18}
 e^{-Tf}D_{b_T}^2e^{Tf}=\left(D_b+T\widehat{c}(\nabla f)\right)^2.
\end{align}
Thus, one has
\begin{align}\label{4.19}
{\rm
 Tr}_s\left[N\exp\left(-tD_{b_T}^2\right) \right] ={\rm
 Tr}_s\left[N\exp\left(-t\left(D_{b }+T\widehat{c}(\nabla f)\right)^2\right) \right]  .
\end{align}

By (\ref{4.19}), one writes
\begin{multline}\label{4.20}
\int_{0}^{d}\left(
 {\rm
 Tr}_s\left[gN\exp\left(-tD_{b_T}^2\right)\right]-
 {a_{-1}\over\sqrt{t}}-{n\over 2}\chi_g(F) \right) {dt\over t}\\
 =2\int_{1}^{ \sqrt{dT}}\left(
 {\rm
 Tr}_s\left[gN\exp\left(-\left({t\over \sqrt{T}}D_{b }+t\sqrt{T}\widehat{c}(\nabla f)\right)^2\right) \right] -
 {\sqrt{T}\over {t}}a_{-1}-{n\over 2}\chi_g(F)\right) {dt\over t}\\
 +2\int_0^{1\over \sqrt{T}}\left(
 {\rm
 Tr}_s\left[gN\exp\left(-\left(tD_{b }+tT\widehat{c}(\nabla f)\right)^2\right) \right] -
 {a_{-1}\over {t}}-{n\over 2}\chi_g(F)\right) {dt\over t} .
\end{multline}

In view of Theorem \ref{t3.5}, we write
\begin{multline}\label{4.21}
\int_0^{1\over \sqrt{T}}\left(
 {\rm
 Tr}_s\left[gN\exp\left(-\left(tD_{b }+tT\widehat{c}(\nabla f)\right)^2\right) \right] -
 {a_{-1}\over {t}}-{n\over 2}\chi_g(F)\right) {dt\over t}\\
 = \int_0^{1\over \sqrt{T}}\left(
 {\rm
 Tr}_s\left[gN\exp\left(-\left(tD_{b }+tT\widehat{c}(\nabla f)\right)^2\right) \right] -
 {{1}\over {t}}\int_{M_g}{\rm Tr}_{F}[g]\int^BL\exp\left(-B_{{(tT)}^2}\right)\right.\\
 \left.-{tT\over 2}\int_{M_g}\theta_g\left(F,b^F\right)\int^B\widehat{df}
 \exp\left(-B_{(tT)^2}\right)-{n\over 2}\chi_g(F)\right) {dt\over t}\\
 + \int_0^{1\over \sqrt{T}}\left({{1}\over {t}}\int_{M_g}{\rm Tr}_{F}[g]\int^BL\exp\left(-B_{{(tT)}^2}\right)
 -{a_{-1}\over t}\right){dt\over t} \\
 +\int_0^{1\over \sqrt{T}}{tT\over 2}\int_{M_g}\theta_g\left(F,b^F\right)\int^B\widehat{df}
 \exp\left(-B_{(tT)^2}\right){dt\over t}.
\end{multline}

By \cite[Definitions 3.6, 3.12 and Theorem 3.18]{BZ1}, one has, as
$T\rightarrow +\infty$,
\begin{multline}\label{4.22}
\int_0^{1\over \sqrt{T}}{tT\over
2}\int_{M_g}\theta_g\left(F,b^F\right)\int^B\widehat{df}
 \exp\left(-B_{(tT)^2}\right){dt\over t}\to\\
 {1\over 2}\int_{M_g}\theta_g\left(F,b^F\right)
 (\nabla f)^*\psi\left(TM_g,\nabla^{TM_g}\right).
\end{multline}

From \cite[(3.54)]{BZ1}, \cite[(3.35)]{SZ} and integration by parts,
we have
\begin{multline}\label{4.23}
\int_0^{1\over \sqrt{T}}\left({{1}\over {t}}\int_{M_g}{\rm
Tr}_{F}[g]\int^BL\exp\left(-B_{{(tT)}^2}\right)
 -{a_{-1}\over t}\right){dt\over t}\\
 =-\sqrt{T}\int_{M_g}{\rm Tr}_{F}[g] \int^B L\exp\left(-B_{ T
}\right)+\sqrt{T}a_{-1}-T\int_{M_g}{\rm
Tr}_{F}[g]f\int^B\exp\left(-B_T\right)
\\ +T\int_{M_g}{\rm Tr}_{F}[g]f\int^B\exp\left(-B_0\right).
\end{multline}
From Theorems \ref{t3.5}, \ref{t3.6}, \cite[Theorem 3.20]{BZ1},
\cite[(7.72) and (7.73)]{BZ1} and the dominate convergence, one
finds that as $T\rightarrow +\infty$,
\begin{multline}\label{4.24}
\int_0^{1\over \sqrt{T}}\left(
 {\rm
 Tr}_s\left[gN\exp\left(-\left(tD_{b }+tT\widehat{c}(\nabla f)\right)^2\right) \right] -
 {{1}\over {t}}\int_{M_g}{\rm Tr}_{F}[g]\int^BL\exp\left(-B_{{(tT)}^2}\right)\right.\\
 \left.-{tT\over 2}\int_{M_g}\theta_g\left(F,b^F\right)\int^B\widehat{df}
 \exp\left(-B_{(tT)^2}\right)-{n\over 2}\chi_g(F)\right) {dt\over t}\\
=\int_0^{1 }\left(
 {\rm
 Tr}_s\left[gN\exp\left(-\left({t\over\sqrt{T}}D_{b }+t\sqrt{T }
 \widehat{c}(\nabla f)\right)^2\right) \right]\right.\\ \left. -
 {\sqrt{T }\over {t}}\int_{M_g}{\rm Tr}_{F}[g]\int^BL\exp\left(-B_{{(t\sqrt{T })}^2}\right)\right.\\
 \left.-{t\sqrt{T }\over 2}\int_{M_g}\theta_g\left(F,b^F\right)\int^B\widehat{df}
 \exp\left(-B_{(t\sqrt{T })^2}\right)-{n\over 2}\chi_g(F)\right) {dt\over
 t}\\
\to\int_0^1\left\{\sum_{j}{\rm Tr}\left[g|_{F|_{M_{g,j}}}\right]\right.\\
 \cdot\left({1\over{1-e^{-2t^2}}}
 \left(\left(1+e^{-2t^2}\right)\sum_{x\in{B\cap{M_{g,j}}}}(-1)^{{\rm ind}_g(x)}{{\rm ind}_g}(x)
 -{\rm dim}M_{g,j}e^{-2t^2}\chi(M_{g,j})\right)\right)\\
 +\sum_{j}{\rm Tr}\left[g|_{F|_{M_{g,j}}}\right]\sum_{k}{\sinh(2t^2)\over{\cosh(2t^2)-\cos(\beta_k)}}\sum_{x\in
 {B\cap{M_{g,j}}}}(-1)^{{\rm ind}_g(x)}n_{-}({\beta_k})(x)\\
 -\sum_{j}{\rm Tr}\left[g|_{F|_{M_{g,j}}}\right]\sum_{k}{1\over
 2}\left({\sinh(2t^2)\over{\cosh(2t^2)-\cos(\beta_k)}}-1\right){\rm
 dim}N^{\beta_k}\chi(M_{g,j})\\
\left.+{1\over{2t^2}}\sum_{j}{\rm
Tr}\left[g|_{F|_{M_{g,j}}}\right]\sum_{x\in{B\cap{M_{g,j}}}}(-1)^{{\rm
ind}_g(x)}\left({\rm dim}M_{g,j}-2{\rm ind}_g(x)\right)-{n\over
2}\chi_g(F)\right\}{dt\over t}\\
={1\over 2}\sum_{j}{\rm
Tr}\left[g|_{F|_{M_{g,j}}}\right]\left\{\sum_{x\in{B\cap{M_{g,j}}}}(-1)^{{\rm
ind}_g(x)}{{\rm ind}_g}(x)-{1\over 2}\sum_{j}\chi(M_{g,j}){\rm
dim}M_{g,j}\right\}\\
\cdot\int_0^1\left({{1+e^{-2t}}\over{1-e^{-2t}}}-{1\over
t}\right){dt\over t}\\
-\sum_{j}{\rm Tr}\left[g|_{F|_{M_{g,j}}}\right]\sum_k\left({1\over
4}{\rm dim}N^{\beta_k}\chi(M_{g,j})-{1\over
2}\sum_{x\in{B\cap{M_{g,j}}}}(-1)^{{\rm
ind}_g(x)}n_{-}(\beta_k)(x)\right)\\
\cdot\int_0^1\left(\sinh(2t)\over{\cosh(2t)-\cos(\beta_k)}\right){dt\over
t}.
\end{multline}

On the other hand, by Theorems \ref{t3.6}, \ref{t3.7} and the
dominate convergence, we have that as $T\to +\infty$,
\begin{multline}\label{4.25}
\int_1^{\sqrt{Td}}\left(
 {\rm
 Tr}_s\left[gN\exp\left(-\left({t\over \sqrt{T}}D_{b }+t\sqrt{T}\widehat{c}(\nabla f)\right)^2\right) \right] -
 {\sqrt{T} \over {t}}a_{-1}-{n\over 2}\chi_g(F)\right) {dt\over t}\\
 = \int_1^{\sqrt{Td}}\left(
 {\rm
 Tr}_s\left[gN\exp\left(-\left({t\over\sqrt{T}}D_{b }+t\sqrt{T}\widehat{c}(\nabla f)\right)^2\right) \right]
 - \widetilde{\chi}'_g(F)
 \right) {dt\over t}\\
 +{1\over 2}  \widetilde{\chi}'_g(F) \log \left(Td\right)+a_{-1}\sqrt{T}\left({1\over\sqrt{Td}}-1\right)
  -{n\over 4}\chi_g(F)\log \left(Td\right)\\
=\int_{1}^{+\infty}\left\{\sum_{j}{\rm Tr}\left[g|_{F|_{M_{g,j}}}\right]\right.\\
 \cdot\left({1\over{1-e^{-2t^2}}}
 \left(\left(1+e^{-2t^2}\right)\sum_{x\in{B\cap{M_{g,j}}}}(-1)^{{\rm ind}_g(x)}{{\rm ind}_g}(x)
 -{\rm dim}M_{g,j}e^{-2t^2}\chi(M_{g,j})\right)\right)\\
 +\sum_{j}{\rm Tr}\left[g|_{F|_{M_{g,j}}}\right]\sum_{k}{\sinh(2t^2)\over{\cosh(2t^2)-\cos(\beta_k)}}\sum_{x\in
 {B\cap{M_{g,j}}}}(-1)^{{\rm ind}_g(x)}n_{-}({\beta_k})(x)\\
 \left.-\sum_{j}{\rm Tr}\left[g|_{F|_{M_{g,j}}}\right]\sum_{k}{1\over
 2}\left({\sinh(2t^2)\over{\cosh(2t^2)-\cos(\beta_k)}}-1\right){\rm
 dim}N^{\beta_k}\chi(M_{g,j})-\widetilde{\chi}'_g(F)\right\}{dt\over t}\\
+{1\over 2}  \widetilde{\chi}'_g(F) \log
\left(Td\right)+a_{-1}\sqrt{T}\left({1\over\sqrt{Td}}-1\right)
  -{n\over 4}\chi_g(F)\log \left(Td\right)+o(1)\\
=\sum_{j}{\rm
Tr}\left[g|_{F|_{M_{g,j}}}\right]\left\{\sum_{x\in{B\cap{M_{g,j}}}}(-1)^{{\rm
ind}_g(x)}{{\rm ind}_g}(x)-{1\over 2}\sum_{j}\chi(M_{g,j}){\rm
dim}M_{g,j}\right\}\\
\cdot\int_{1}^{+\infty}{e^{-2t}\over{1-e^{-2t}}}{dt\over t}\\
-\sum_{j}{\rm Tr}\left[g|_{F|_{M_{g,j}}}\right]\sum_k\left({1\over
4}{\rm dim}N^{\beta_k}\chi(M_{g,j})-{1\over
2}\sum_{x\in{B\cap{M_{g,j}}}}(-1)^{{\rm
ind}_g(x)}n_{-}(\beta_k)(x)\right)\\
\cdot\int_1^{+\infty}\left({\sinh(2t)\over{\cosh(2t)-\cos(\beta_k)}}-1\right){dt\over
t}\\ +{1\over 2}\left(\widetilde{\chi}'_g(F)-{n\over
2}\chi_g(F)\right)\log(Td)+{a_{-1}\over\sqrt{d}}-\sqrt{T}a_{-1}+o(1).
\end{multline}

Combining (\ref{3.4}), (\ref{4.14}) and (\ref{4.20})-(\ref{4.25}),
one deduces, by setting $T\to +\infty$, that
\begin{multline}\label{4.26}
\log\left({P_\infty^{\det H}\left(b^{\rm RS}_{{\rm
det}(H^*(M,F),G)}\right)\over b_{{\rm det}(H^*(W^u,F),G)}^{{\cal
M},-X}}(g) \right)\simeq\\
-2{\rm Tr}_{s}^{B_g}[f]T+\left(\widetilde{\chi}'_{g}(F)-{n\over
2}\chi_{g}(F)\right)\log T-\left(\widetilde{\chi}'_{g}(F)-{n\over
2}\chi_{g}(F)\right)\log\pi\\
-\int_{M_g}\theta_g\left(F,b^F\right)
 (\nabla f)^*\psi\left(TM_g,\nabla^{TM_g}\right)\\
+2\sqrt{T}\int_{M_g}{\rm Tr}_{F}[g] \int^B L\exp\left(-B_{ T
}\right)-2\sqrt{T}a_{-1}+2T\int_{M_g}{\rm
Tr}_{F}[g]f\int^B\exp\left(-B_T\right)
\\ -2T\int_{M_g}{\rm Tr}_{F}[g]f\int^B\exp\left(-B_0\right)\\
-\sum_{j}{\rm
Tr}\left[g|_{F|_{M_{g,j}}}\right]\left\{\sum_{x\in{B\cap{M_{g,j}}}}(-1)^{{\rm
ind}_g(x)}{{\rm ind}_g}(x)-{1\over 2}\sum_{j}\chi(M_{g,j}){\rm
dim}M_{g,j}\right\}\\
\cdot\left(\int_0^1\left({{1+e^{-2t}}\over{1-e^{-2t}}}-{1\over
t}\right){dt\over t}+\int_{1}^{+\infty}{2e^{-2t}\over{1-e^{-2t}}}{dt\over t}\right)\\
+2\sum_{j}{\rm Tr}\left[g|_{F|_{M_{g,j}}}\right]\sum_k\left({1\over
4}{\rm dim}N^{\beta_k}\chi(M_{g,j})-{1\over
2}\sum_{x\in{B\cap{M_{g,j}}}}(-1)^{{\rm
ind}_g(x)}n_{-}(\beta_k)(x)\right)\\
\cdot\left(\int_0^1\left(\sinh(2t)\over{\cosh(2t)-\cos(\beta_k)}\right){dt\over
t}+\int_1^{+\infty}\left({\sinh(2t)\over{\cosh(2t)-\cos(\beta_k)}}-1\right){dt\over
t}\right)\\
 -\left(\widetilde{\chi}'_g(F)-{n\over
2}\chi_g(F)\right)\log(Td)-2{a_{-1}\over\sqrt{d}}+2\sqrt{T}a_{-1}\\
+2T\int_{M_g}{\rm
Tr}_{F}[g]fe\left(TM_g,\nabla^{TM_g}\right)+{2a_{-1}\over\sqrt{d}}-\left(\Gamma'(1)-\log
d\right)\left(\widetilde{\chi}'_{g}(F)-{n\over
2}\chi_{g}(F)\right)+o(1).\\
\end{multline}

By \cite[Theorem 3.20]{BZ1} and \cite[(7.72)]{BZ1}, one has
\begin{multline}\label{4.27}
\lim_{T\to +\infty}\left(2T\int_{M_g}{\rm
Tr}_{F}[g]f\int^{B}\exp(-B_T)-2T{\rm Tr}_{s}^{B_g}[f]\right)\\
=-\sum_{j}{\rm
Tr}\left[g|_{F|_{M_{g,j}}}\right]\left\{\sum_{x\in{B\cap{M_{g,j}}}}(-1)^{{\rm
ind}_g(x)}{{\rm ind}_g}(x)-{1\over 2}\sum_{j}\chi(M_{g,j}){\rm
dim}M_{g,j}\right\},
\end{multline}
\begin{multline}\label{4.28}
\lim_{T\to +\infty}2\sqrt{T}\int_{M_g}{\rm
Tr}_{F}[g]\int^{B}L\exp(-B_T)\\
=2\sum_{j}{\rm
Tr}\left[g|_{F|_{M_{g,j}}}\right]\left\{\sum_{x\in{B\cap{M_{g,j}}}}(-1)^{{\rm
ind}_g(x)}{{\rm ind}_g}(x)-{1\over 2}\sum_{j}\chi(M_{g,j}){\rm
dim}M_{g,j}\right\}.
\end{multline}

On the other hand, by \cite[(7.93)]{BZ1} and \cite[(5.55)]{BZ2}, one
has
\begin{align}\label{4.29}
\int_0^1\left({1+e^{-2t}\over 1-e^{-2t}}-{1\over
 t}\right){dt\over t}+\int_1^{+\infty} {2\, e^{-2t}\over 1-e^{-2t}} {dt\over
 t}  =1-\log\pi-\Gamma'(1),
\end{align}
\begin{multline}\label{4.30}
\int_0^1\left(\sinh(2t)\over{\cosh(2t)-\cos(\beta_k)}\right){dt\over
t}+\int_1^{+\infty}\left({\sinh(2t)\over{\cosh(2t)-\cos(\beta_k)}}-1\right){dt\over
t}\\
=-\log(\pi)-{1\over
2}\left({\Gamma'\over\Gamma}\left(\beta_k\over{2\pi}\right)+{\Gamma'\over\Gamma}\left(1-{\beta_k\over{2\pi}}\right)\right).\\
\end{multline}
Also, by \cite[(5.64)]{BZ2}, if $x\in{B\cap{M_g}}$,
\begin{align}\label{4.31}
{{\rm dim}N^{\beta_k}\over 4}-{{n_{-}(\beta_k)(x)}\over 2}={1\over
4}\left[n_{+}(\beta_k)(x)-n_{-}(\beta_k)(x)\right].
\end{align}
From (\ref{4.26})-(\ref{4.31}), we get (\ref{2.54}), which completes
the proof of Theorem \ref{t2.10}.

\section{Proofs of the intermediary Theorems} \label{s5}
\setcounter{equation}{0}

The purpose of this section is to give a sketch of the proofs of the
intermediary Theorems. Since the methods of the proofs of these
theorems are essentially the same as the corresponding theorem in
\cite{SZ}, so we will refer to \cite{SZ} for related definitions and
notations directly when there will be no confusion, such as
$B_{b,g}$, $A_{b,t,T}$, $A_{g,t,T}$, $C_{t,T}$, $\cdots$.

\subsection{Proof of Theorem \ref{t3.1}}
\label{s5.1}

From Theorem \ref{t2.5} and \cite[(4.44)]{SZ} which in our situation
we also have that $P_{\infty, T}$ commutate with $g\in G$, one finds
\begin{align}\label{5.*1}
{P_T^{[0,1],\det H}\left(b^{\rm RS}_{\det (H^*(\Omega^*_{[0,1],T}
(M,F)),G)}\right)\over b_{{\rm det}(H^*(W^u,F),G)}^{{\cal M},-X}
}(g)=\prod_{i=0}^n
\det\left(\left.P_{\infty,T}^\#P_{\infty,T}\right|_{\Omega^i_{[0,1],T}(M,F)}\right)^{(-1)^{i+1}}(g).
\end{align}

From \cite[Propositions 4.4 and 4.5]{SZ}, one deduces that as
$T\rightarrow+\infty$,
\begin{multline}\label{5.*2}
\det\left(\left.P_{\infty,T}^\#P_{\infty,T}\right|_{\Omega^i_{[0,1],T}(M,F)}\right)(g)
\\
=\det\left(\left.e_Te_T^\#P_{\infty,T}^\#P_{\infty,T}\right|_{\Omega^i_{[0,1],T}(M,F)}\right)(g)\cdot
{\det}^{-1}\left(\left.e_Te_T^\#
\right|_{\Omega^i_{[0,1],T}(M,F)}\right)(g)\\
=\det\left(\left.
\left(P_{\infty,T}e_T\right)^\#P_{\infty,T}e_T\right|_{C^i(W^u,F)}\right)(g)\cdot
{\det}^{-1}\left(\left.e_T^\#e_T \right|_{C^i(W^u,F)}\right)(g)\\
=\det\left(\left.\left(1+O\left(e^{-cT}\right)\right)^\#\left({\pi\over
T}\right)^{N-n/2}e^{2T{\cal F}}\left(1+O\left(e^{-cT}\right)\right)
 \right|_{C^i(W^u,F)}\right)(g)\\ \cdot
{\det}^{-1}\left(\left.\left(1+O\left(e^{-cT}\right)\right)\right|_{C^i(W^u,F)}\right)(g).
\end{multline}

From (\ref{5.*1}) and (\ref{5.*2}), one gets (\ref{3.4})
immediately.

The proof of Theorem \ref{t3.1} is completed. \ \ Q.E.D.

\subsection{Proof of Theorem \ref{t3.2}}
\label{s5.2}

The proof of Theorem \ref{t3.2} is the same as the proof of
\cite[Theorem 3.4]{SZ} given in \cite[Section 5]{SZ}.

\subsection{Proof of Theorem \ref{t3.3}}
\label{s5.3}

Recall that the operator $e_T:C^*(W^u,F)\to \Omega_{[0,1],T}^*
(M,F)$ has been defined in \cite[(4.38)]{SZ}, and in the current
case, we also have that $e_T$ commute with $G$. So by
\cite[Proposition 4.4]{SZ}, we have that for $T\geq 0$ large enough,
$e_T : C^*(W^u,F)\to \Omega_{[0,1],T}^* (M,F)$ is an identification
of $G$-spaces. So (\ref{3.7}) follows. Also (\ref{3.8}) was already
proved in \cite[Theorem 3.5]{SZ}.

\subsection{Proof of Theorem \ref{t3.4}}
\label{s5.4}

In this section, we provide a proof of Theorem \ref{t3.4}, which
computes the asymptotic of ${\rm Tr}_{s}[gN\exp(-tD_{b_T}^{2})]$ for
fixed $T\geq 0$ as $t\to 0$. The method is the essentially same as
it in \cite{SZ}.

By \cite[(6.4)]{SZ}, we have
\begin{multline}\label{5.1}
e^{-tD_b^2}=e^{-tD_g^2}+\sum_{k=1}^{n}(-1)^kt^k\int_{\Delta_k}e^{-t_1tD_g^2}B_{b,g}
e^{-t_2tD_g^2}\cdots B_{b,g}e^{-t_{k+1}tD_g^2}dt_1\cdots dt_k\\
+ (-1)^{n+1}t^{n+1}\int_{\Delta_{n+1}}e^{-t_1tD_g^2}B_{b,g}
e^{-t_2tD_g^2}\cdots B_{b,g}e^{-t_{n+2}tD_b^2}dt_1\cdots dt_{n+1},
\end{multline}
where $\Delta_k$, $1\leq k\leq n+1$, is the $k$-simplex defined by
$t_1+\cdots+t_{k+1}=1$, $t_1\geq 0$, $\cdots, $ $t_{k+1}\geq 0$.
Also, by the same proof of \cite[Proposition 6.1]{SZ}, we have the
following result.
\begin{prop}\label{t5.1}
As $t\rightarrow 0^+$, one has
\begin{align}\label{5.2}
t^{n+1}\int_{\Delta_{n+1}}{\rm Tr}_s\left[gNe^{-t_1tD_g^2}B_{b,g}
e^{-t_2tD_g^2}\cdots B_{b,g}e^{-t_{n+2}tD_b^2}\right]dt_1\cdots
dt_{n+1}\rightarrow 0.
\end{align}
\end{prop}
By \cite[(6.22) and (6.23)]{SZ}, we have that for any $ 1< k\leq n$,
$(t_1,\cdots,t_{k+1})\in\Delta_k$,
\begin{align}\label{5.3}
\lim_{t\rightarrow 0^+}t^k {\rm Tr}_s\left[gNe^{-t_1tD_g^2}B_{b,g}
e^{-t_2tD_g^2}\cdots B_{b,g}e^{-t_{k+1}tD_g^2}\right]=0,
\end{align}
while for $k=1, 0\leq t_1\leq 1$,
\begin{multline}\label{5.4}
\lim_{t\rightarrow 0^+}t {\rm Tr}_s\left[gNe^{-t_1tD_g^2}B_{b,g}
 e^{-(1-t_1)tD_g^2}\right]=\lim_{t\rightarrow 0^+}t {\rm Tr}_s\left[gNB_{b,g}
 e^{- tD_g^2}\right]\\
={1\over 2}\int_{M_g}\int^B {\rm Tr}\left[g\left( \sum_{i,\, j=1}^n
   e_i \wedge\widehat{e_j}\left(\nabla_{e_i}^u\omega^F\left(e_j\right)\right)+{1\over
  2}\left[
   \omega^F ,\widehat{ \omega^F_g}-\widehat{ \omega^F }
   \right]\right)\right]\\
   \cdot L\exp\left(-{\dot{R}^{TM_g}\over 2}\right).
\end{multline}
So by \cite[(2.13)]{BZ2}, and proceed as in
\cite[(6.26)-(6.28)]{SZ}, we have
\begin{align}\label{5.5}
\lim_{t\rightarrow 0^+}t {\rm Tr}_s\left[gNe^{-t_1tD_g^2}B_{b,g}
 e^{-(1-t_1)tD_g^2}\right]=0.
\end{align}

From (\ref{5.1}), (\ref{5.2}), (\ref{5.3}), (\ref{5.5}) and
\cite[Theorem 5.9]{BZ2}, one gets (\ref{3.9}).

The proof of Theorem \ref{t3.4} is completed. \ \ Q.E.D.

\subsection{Proof of Theorem \ref{t3.5}}
\label{s5.5}

In order to prove (\ref{3.10}), one need only to prove that under
the conditions of Theorem \ref{t3.5}, there exists constant $C''>0$
such that
\begin{multline}\label{5.6}
 \left|  {\rm Tr}_s\left[gN\exp\left(-\left(tD_{b }+T\widehat{c}
 (\nabla f)\right)^2\right)\right]- {\rm Tr}_s\left[gN\exp\left(-\left(tD_{g}+T\widehat{c}
 (\nabla f)\right)^2\right)\right]\right.\\
\left.  -{T\over
2}\int_{M_g}\left(\theta_{g}\left(F,b^F\right)-\theta_{g}\left(F,g^F\right)\right)
\int^B\widehat{df}\exp\left(-B_{T^2}\right) \right|\leq C''t.
\end{multline}

By \cite[(7.8)]{SZ}, we have
\begin{multline}\label{5.7}
 e^{-A_{b,t,T}^2}=e^{-A_{g,t,T}^2}\\
 +\sum_{k=1}^{n}(-1)^k\int_{\Delta_k}e^{-t_1A_{g,t,T}^2}C_{t,T}
e^{-t_2A_{g,t,T}^2}\cdots C_{t,T}e^{-t_{k+1}A_{g,t,T}^2}dt_1\cdots dt_k\\
+ (-1)^{n+1} \int_{\Delta_{n+1}}e^{-t_1A_{g,t,T}^2}C_{t,T}
e^{-t_2A_{g,t,T}^2}\cdots C_{t,T}e^{-t_{n+2}A_{b,t,T}^2}dt_1\cdots
dt_{n+1}.
\end{multline}
By the same proof of \cite[(7.21)]{SZ}, we have that there exists
$C_1>0$ such that for any $t>0$ small enough and $T\in[0,{1\over
t}]$,
\begin{align}\label{5.8}
\left|\int_{\Delta_{n+1}}{\rm Tr}_s\left[g
Ne^{-t_1A_{g,t,T}^2}C_{t,T} e^{-t_2A_{g,t,T}^2}\cdots
C_{t,T}e^{-t_{n+2}A_{b,t,T}^2}\right]dt_1\cdots dt_{n+1}\right| \leq
{C_1 t }   .
\end{align}
Also by the same proof of \cite[(7.23)]{SZ}, we have that there
exists $C_2>0$, $0<d<1$ such that for any $1<k\leq n$, $0<t\leq d$,
$T\geq 0$ with $tT\leq 1$,
\begin{align}\label{5.9}
\left|\int_{\Delta_k}{\rm Tr}_s\left[gNe^{-t_1A_{g,t,T}^2}C_{t,T}
e^{-t_2A_{g,t,T}^2}\cdots
C_{t,T}e^{-t_{k+1}A_{g,t,T}^2}\right]dt_1\cdots dt_k\right|\leq
C_2t,
\end{align}
while for $k=1$ one has for any $0<t\leq d$, $T\geq 0$ with $tT\leq
1$ and $0\leq t_1\leq 1$, by \cite[Proposition 9.3]{BZ2}, we have
\begin{multline}\label{5.10}
\left| {\rm Tr}_s\left[gNe^{-t_1A_{g,t,T}^2}C_{t,T}
e^{-\left(1-t_1\right)A_{g,t,T}^2} \right] - T\int_{M_g}\int^B{\rm
Tr}\left[g\omega^F(\nabla
f)\right]L\exp\left(-B_{T^2}\right)\right|\\ \leq C_2t.
\end{multline}

Now similar as \cite[(7.25)]{SZ}, we have
\begin{multline}\label{5.11}
\int_{M_g}\int^B{\rm Tr}\left[g\omega^F(\nabla
f)\right]L\exp\left(-B_{T^2}\right)\\= {1\over 2}\int_{M_g}
\left(\theta_{g}\left(F,g^F\right)-\theta_{g}\left(F,b^F\right)\right)
\int^B \widehat{ \nabla f} \exp\left(-B_{T^2}\right).
\end{multline}

From (\ref{5.7})-(\ref{5.11}), we get (\ref{5.6}), which completes
the proof of Theorem \ref{t3.5}. \ \ Q.E.D.

\subsection{Proof of Theorem \ref{t3.6}}
\label{s5.6}

In order to prove Theorem \ref{t3.6}, we need only to prove that for
any $T>0$,
\begin{align}\label{5.12}
\lim_{t\rightarrow 0^+}\left({\rm
Tr}_s\left[gN\exp\left({-A_{b,t,{T\over t}}^2}\right)\right] -{\rm
Tr}_s\left[gN\exp\left({-A_{g,t,{T\over
t}}^2}\right)\right]\right)=0.
\end{align}

By \cite[(8.2) and (8.4)]{SZ}, there exists $0<C_0\leq 1$, such that
when $0<t\leq C_0$, one has the absolute convergent expansion
formula
\begin{multline}\label{5.13}
e^{-A_{b,t,{T\over t}}^2}-e^{-A_{g,t,{T\over t}}^2}\\
 =\sum_{k=1}^{+\infty}(-1)^k\int_{\Delta_k}e^{-t_1A_{g,t,{T\over t}}^2}C_{t,{T\over t}}
e^{-t_2A_{g,t,{T\over t}}^2}\cdots C_{t,{T\over
t}}e^{-t_{k+1}A_{g,t,{T\over t}}^2}dt_1\cdots dt_k ,
\end{multline}
and that
\begin{align}\label{5.14}
\sum_{k=n}^{+\infty}(-1)^k\int_{\Delta_k}g e^{-t_1A_{g,t,{T\over
t}}^2}C_{t,{T\over t}} e^{-t_2A_{g,t,{T\over t}}^2}\cdots
C_{t,{T\over t}}e^{-t_{k+1}A_{g,t,{T\over t}}^2}dt_1\cdots dt_k
\end{align}
is uniformly absolute convergent for $0<t\leq C_0$.

Proceed as in \cite[Section 8]{SZ}, one has that for any
$(t_1,\cdots, t_{k+1})\in\Delta_k\setminus \{t_1\cdots t_{k+1}=0\}$,
\begin{multline}\label{5.15}
\left|{\rm Tr}_s\left[g N  e^{-t_1A_{g,t,{T\over t}}^2}C_{t,{T\over
t}} e^{-t_2A_{g,t,{T\over t}}^2}\cdots C_{t,{T\over
t}}e^{-t_{k+1}A_{g,t,{T\over t}}^2}\right]\right| \\
\leq C_3t^k  \left(t_1\cdots t_{k}\right)^{-{1\over 2}} {\rm
Tr}\left[e^{-{A_{g,t,{T\over t}}^2\over 2}}\right] \left\| \psi
e^{-{ t_{k+1}\over 2}A_{g,t,{T\over t}}^2} \right\|
\end{multline}
for some positive constant $C_3 >0$.

Also, by \cite[(8.4)]{SZ}, (\ref{5.15}) and the same assumption in
\cite{SZ} that $t_{k+1}\geq{1\over{k+1}}$, one gets
\begin{multline}\label{5.16}
\left|\int_{\Delta_k}{\rm Tr}_s\left[g N  e^{-t_1A_{g,t,{T\over
t}}^2}C_{t,{T\over t}} e^{-t_2A_{g,t,{T\over t}}^2}\cdots
C_{t,{T\over t}}e^{-t_{k+1}A_{g,t,{T\over t}}^2}\right]dt_1\cdots dt_k\right| \\
 \leq C_4  t^{k-n}    \left\| \psi
e^{-{ 1\over 2(k+1)}A_{g,t,{T\over t}}^2} \right\|
\end{multline}
for some constant $C_4 >0$.

From (\ref{5.13}), (\ref{5.14}), (\ref{5.16}), \cite[(8.9) and
(8.10)]{SZ} and the dominate convergence, we get (\ref{5.12}), which
completes the proof of Theorem \ref{t3.6}. \ \ Q.E.D.

\subsection{Proof of Theorem \ref{t3.7}}
\label{s5.7}

In order to prove Theorem \ref{t3.7}, we need only to prove that
there exist $c>0$, $C>0$, $0<C_0\leq 1$ such that for any  $0<t\leq
C_0$, $ T\geq 1$,
\begin{align}\label{5.17}
 \left| {\rm
Tr}_s\left[gN\exp\left({-A_{b,t,{T\over t}}^2}\right)\right] -{\rm
Tr}_s\left[gN\exp\left({-A_{g,t,{T\over
t}}^2}\right)\right]\right|\leq c\exp(-CT).
\end{align}

First of all, one can choose $C_0>0$ small enough so that  for any
$0<t\leq C_0$, $T>0$, by (\ref{5.13}), we have the  absolute
convergent expansion formula
\begin{multline}\label{5.18}
e^{-A_{b,t,{T\over t}}^2}-e^{-A_{g,t,{T\over t}}^2}\\
 =\sum_{k=1}^{+\infty}(-1)^k\int_{\Delta_k}e^{-t_1A_{g,t,{T\over t}}^2}C_{t,{T\over t}}
e^{-t_2A_{g,t,{T\over t}}^2}\cdots C_{t,{T\over
t}}e^{-t_{k+1}A_{g,t,{T\over t}}^2}dt_1\cdots dt_k ,
\end{multline}
from which one has
\begin{multline}\label{5.19}
{\rm Tr}_s\left[gN\exp\left({-A_{b,t,{T\over t}}^2}\right)\right]
-{\rm Tr}_s\left[gN\exp\left({-A_{g,t,{T\over
t}}^2}\right)\right]\\
 =\sum_{k=1}^{+\infty}(-1)^k\int_{\Delta_k}{\rm Tr}_s\left[g N
 e^{-t_1A_{g,t,{T\over t}}^2}C_{t,{T\over t}}
e^{-t_2A_{g,t,{T\over t}}^2}\cdots C_{t,{T\over
t}}e^{-t_{k+1}A_{g,t,{T\over t}}^2}\right]dt_1\cdots dt_k .
\end{multline}
Thus, in order to prove (\ref{5.17}), we need only to prove
\begin{multline}\label{5.20}
\sum_{k=1}^{+\infty}\left|\int_{\Delta_k}{\rm Tr}_s\left[g N
 e^{-t_1A_{g,t,{T\over t}}^2}C_{t,{T\over t}}
e^{-t_2A_{g,t,{T\over t}}^2}\cdots C_{t,{T\over
t}}e^{-t_{k+1}A_{g,t,{T\over t}}^2}\right]dt_1\cdots dt_k\right|\\
=\sum_{k=1}^{+\infty}\left|\int_{\Delta_k}{\rm Tr}_s\left[g N
 e^{-\left(t_1+t_{k+1}\right)A_{g,t,{T\over t}}^2}C_{t,{T\over t}}
e^{-t_2A_{g,t,{T\over t}}^2}\cdots C_{t,{T\over t}}\right]dt_1\cdots
dt_k\right|\\ \leq c\exp(-CT) .
\end{multline}

By \cite[(8.6)]{SZ}, we have for any $t>0$, $T\geq 1$,
$(t_1,\cdots,t_{k+1})\in \Delta_k\setminus\{t_1\cdots t_{k+1}=0\}$,
\begin{multline}\label{5.21}
{\rm Tr}_s\left[g N
 e^{-\left(t_1+t_{k+1}\right)A_{g,t,{T\over t}}^2}C_{t,{T\over t}}
e^{-t_2A_{g,t,{T\over t}}^2}\cdots C_{t,{T\over
t}}\right] \\
={\rm Tr}_s\left[g N
 \psi e^{-\left(t_1+t_{k+1}\right)A_{g,t,{T\over t}}^2}  C_{t,{T\over t}}
\psi e^{-t_2A_{g,t,{T\over t}}^2}  C_{t,{T\over t}}\cdots\psi
e^{-t_kA_{g,t,{T\over t}}^2}  C_{t,{T\over t}}\right] .
\end{multline}

From (\ref{5.21}), \cite[(9.18) and (9.19)]{SZ}, one sees that there
exists $C_5>0$, $C_6>0$ and $C_7>0$ such that for any $k\geq 1$,
\begin{multline}\label{5.22}
\left|\int_{\Delta_k}{\rm Tr}_s\left[g N
 e^{-t_1A_{g,t,{T\over t}}^2}C_{t,{T\over t}}
e^{-t_2A_{g,t,{T\over t}}^2}\cdots C_{t,{T\over
t}}e^{-t_{k+1}A_{g,t,{T\over t}}^2}\right]dt_1\cdots dt_k\right|\\
\leq C_{5} \left(C_6t\right)^k{T^{n\over 2}\over t^n}
  \exp\left(-{C_7T\over 4}\right),
\end{multline}
from which one sees that there exists $0<c_1\leq 1$, $C_8>0$,
$C_9>0$ such that for any $0<t\leq c_1$ and $T\geq 1$, one has
\begin{multline}\label{5.23}
\left|\sum_{k=n}^{+\infty}\int_{\Delta_k}{\rm Tr}_s\left[g N
 e^{-t_1A_{g,t,{T\over t}}^2}C_{t,{T\over t}}
e^{-t_2A_{g,t,{T\over t}}^2}\cdots C_{t,{T\over
t}}e^{-t_{k+1}A_{g,t,{T\over t}}^2}\right]dt_1\cdots dt_k\right|\\
\leq C_{8}
  \exp\left(-C_{9}T\right).
\end{multline}

On the other hand, for any $1\leq k<n$, by proceeding as in
(\ref{5.16}), one has that for any $0<t\leq c_1$, $T\geq 1$,
\begin{multline}\label{5.24}
\left|\int_{\Delta_k}{\rm Tr}_s\left[g N  e^{-t_1A_{g,t,{T\over
t}}^2}C_{t,{T\over t}} e^{-t_2A_{g,t,{T\over t}}^2}\cdots
C_{t,{T\over t}}e^{-t_{k+1}A_{g,t,{T\over t}}^2}\right]dt_1\cdots dt_k\right| \\
 \leq C_{10}  t^{k-n}    \left\| \psi
e^{-{ 1\over 2(k+1)}A_{g,t,{T\over t}}^2} \right\|
\end{multline}
for some constant $C_{10}>0$.

From (\ref{5.24}) and \cite[(9.23)]{SZ}, one sees immediately that
there exists $C_{11}>0$, $C_{12}>0$ such that for any $1\leq k\leq
n-1$, $0<t\leq c_1$ and $T\geq 1$, one has
\begin{multline}\label{5.25}
\left|\int_{\Delta_k}{\rm Tr}_s\left[g N  e^{-t_1A_{g,t,{T\over
t}}^2}C_{t,{T\over t}} e^{-t_2A_{g,t,{T\over t}}^2}\cdots
C_{t,{T\over t}}e^{-t_{k+1}A_{g,t,{T\over t}}^2}\right]dt_1\cdots dt_k\right| \\
 \leq C_{11}
e^{-C_{12}T } .
\end{multline}

From (\ref{5.19}), (\ref{5.23}) and (\ref{5.25}), one gets
(\ref{5.17}).

The proof of Theorem \ref{t3.7} is completed. \ \ Q.E.D.


\begin {thebibliography}{15}



\bibitem[BGS]{BGS} J.-M. Bismut, H. Gillet and C. Soul\'e,
Analytic torsions and holomorphic determinant line bundles I. {\it
Commun. Math. Phys.} 115 (1988), 49-78.


\bibitem[BGV]{BGV} N. Berline, E. Getzler and M. Vergne, {\it Heat Kernels and Dirac
Operators.} Springer, Berline-Heidelberg-New York, 1992.







\bibitem[BZ1]{BZ1} J.-M. Bismut and W. Zhang, {\it An Extension of a
Theorem by Cheeger and M\"uller}. {\it Ast\'erisque} Tom. 205,
Paris, (1992).

 \bibitem[BZ2]{BZ2} J.-M. Bismut and W. Zhang,
Milnor and Ray-Singer metrics on the equivariant determinant of a
flat vector bundle. {\it Geom.  Funct. Anal.} 4 (1994), 136-212.




 \bibitem [BH1]{BH1} D. Burghelea and  S. Haller,  Torsion, as function on the space of representations.
  {\it Preprint,} math.DG/0507587.

 \bibitem [BH2]{BH2} D. Burghelea and  S. Haller, Complex valued
 Ray-Singer torsion. {\it Preprint,} math.DG/0604484.


\bibitem [BH3]{BH3} D. Burghelea and  S. Haller, Complex valued
 Ray-Singer torsion II.
 {\it Preprint,} math.DG/0610875.



 \bibitem[C]{C} J. Cheeger, Analytic torsion and the heat equation.
 {\it Ann. of Math.} 109 (1979), 259-332.



\bibitem[FT]{FT} M. Farber and V. Turaev, Poincar\'e-Reidemeister
metric, Euler structures and torsion. {\it J. Reine Angew. Math.}
520 (2000), 195-225.





 \bibitem [HS]{HS} B. Helffer and J. Sj\"ostrand, Puis multiples
 en m\'ecanique semi-classique IV: Etude du complexe de Witten.
 {\it Comm. PDE} 10 (1985), 245-340.

\bibitem [KM]{KM} F. F. Knudson and D. Mumford, The projectivity
of the moduli space of stable curves I: Preliminaries on ``det'' and
``div''. {\it Math. Scand.} 39 (1976), 19-55.

\bibitem [L]{La} F. Laudenbach, On the Thom-Smale complex.
Appendix in [BZ1].





\bibitem[MQ]{MQ}   V. Mathai and D. Quillen, Superconnections,
Thom classes, and equivariant differential forms. {\it Topology} 25
(1986), 85-110.


 \bibitem[Mi]{Mi} J. Milnor, Whitehead torsion. {\it  Bull.
 Amer.  Math. Soc.} 72 (1966), 358-426.

 \bibitem[Mu1]{Mu1}  W. M\"uller,  Analytic torsion  and the  R-torsion
 of Riemannian manifolds. {\it Adv. in Math.} 28 (1978), 233-305.

  \bibitem[Mu2]{Mu2}  W. M\"uller,  Analytic torsion  and the  R-torsion
  for unimodular representations. {\it J. Amer. Math. Soc.} 6
  (1993), 721-753.


  \bibitem [Q]{Q}  D. Quillen, Determinants of Cauchy-Riemann
  operators over a Riemann surface. {\it Funct. Anal. Appl.} 14
  (1985), 31-34.



\bibitem [RS]{RS} D. B. Ray and I. M. Singer, $R$-torsion and the
Laplacian on Riemannian manifolds. {\it Adv. in Math.} 7 (1971),
145-210.



\bibitem [S]{S} M. A. Shubin, {\it Pseudodifferential Operators
and Spectral Operator}. Springer-Verlag, Berlin, 2001.



\bibitem [Sm]{Sm} S. Smale, On gradient dynamical systems. {\it Ann.
of Math.} 74 (1961), 199-206.

\bibitem [SZ]{SZ} G. Su and W. Zhang,  A Cheeger-M\"uller theorem
for symmetric bilinear torsions. {\it Preprint}, math.DG/0610577.

\bibitem[T]{T} V. Turaev, Euler structures, nonsingular vector
fields, and Reidemeister-type torsion. {\it Math. USSR-Izv.} 34
(1990), 627-662.

 \bibitem[W]{W}  E. Witten, Supersymmetry and Morse theory. {\it
 J. Diff. Geom.} 17 (1982), 661-692.


 \bibitem[Z]{Z} W. Zhang, {\it Lectures on Chern-Weil Theory and Witten Deformations,}
      Nankai Tracts in Mathematics, Vol. 4.         World Scientific, Singapore, 2001.




\end{thebibliography}
\end{document}